\documentclass[11pt, headinclude,footinclude, BCOR5mm]{scrartcl}
\usepackage[
nochapters, 
beramono, 
eulermath,
pdfspacing, 
dottedtoc 
]{classicthesis} 

\usepackage{arsclassica} 

\usepackage[T1]{fontenc} 

\usepackage[utf8]{inputenc} 

\usepackage{graphicx} 
\graphicspath{{Figures/}} 

\usepackage{enumitem} 

\usepackage{lipsum} 

\usepackage{subfig} 

\usepackage{amsmath,amssymb,amsthm} 

\usepackage{varioref} 

\theoremstyle{definition} 

\theoremstyle{plain} 

\theoremstyle{remark} 

\hypersetup{
colorlinks=true, breaklinks=true, bookmarks=true,bookmarksnumbered,
urlcolor=webbrown, linkcolor=RoyalBlue, citecolor=webgreen, 
pdftitle={}, 
pdfauthor={\textcopyright}, 
pdfsubject={}, 
pdfkeywords={}, 
pdfcreator={pdfLaTeX}, 
pdfproducer={LaTeX with hyperref and ClassicThesis} 
}
\usepackage{titlesec}
\titlelabel{\thetitle.\quad}
\usepackage[utf8]{inputenc}
\usepackage[english]{babel}
\usepackage{enumitem}
\usepackage{amsmath}
\usepackage{amssymb}
\usepackage{xcolor}
\usepackage{lipsum}
\usepackage{longtable}
\usepackage{amsthm}
\usepackage{mathtools}
\usepackage{amsfonts}
\usepackage{graphicx}
\usepackage{float}
\usepackage{dsfont}

\usepackage{algorithm}
\usepackage{algorithmic}
\makeatletter
\def\listofalgorithms{\@starttoc{loa}\listalgorithmname}
\def\l@algorithm{\@tocline{0}{3pt plus2pt}{0pt}{1.9em}{}}
\renewcommand{\ALG@name}{Algorithm}
\renewcommand{\listalgorithmname}{List of \ALG@name s}
\numberwithin{algorithm}{section}

\usepackage[margin=2cm]{geometry}
\theoremstyle{plain}
\numberwithin{equation}{section}
\newtheorem{dfn}{Definition}[section]

\newtheorem{thm}[dfn]{Theorem}

\newtheorem{alg}[dfn]{Algorithm}

\theoremstyle{definition}

\newtheorem{exm}[dfn]{Example}

\newcommand{\norm}[1]{\left\lVert#1\right\rVert}

\usepackage{color, colortbl}
\definecolor{Gray}{gray}{0.9}
\newcolumntype{H}{>{\setbox0=\hbox\bgroup}c<{\egroup}@{}}

\usepackage{etoolbox}
\patchcmd{\thebibliography}
  {\settowidth}
  {\setlength{\itemsep}{0pt plus 0.1pt}\settowidth}
  {}{}
\apptocmd{\thebibliography}
  {\small}
  {}{}

\makeatletter
\DeclareOldFontCommand{\rm}{\normalfont\rmfamily}{\mathrm}
\DeclareOldFontCommand{\sf}{\normalfont\sffamily}{\mathsf}
\DeclareOldFontCommand{\tt}{\normalfont\ttfamily}{\mathtt}
\DeclareOldFontCommand{\bf}{\normalfont\bfseries}{\mathbf}
\DeclareOldFontCommand{\it}{\normalfont\itshape}{\mathit}
\DeclareOldFontCommand{\sl}{\normalfont\slshape}{\@nomath\sl}
\DeclareOldFontCommand{\sc}{\normalfont\scshape}{\@nomath\sc}
\makeatother
\makeatletter
\def\@seccntformat#1{\@ifundefined{#1@cntformat}%
   {\csname the#1\endcsname\space}
   {\csname #1@cntformat\endcsname}}
\newcommand\section@cntformat{\thesection.\space}       
\newcommand\subsection@cntformat{\thesubsection.\space} 
\makeatother
\allowdisplaybreaks
\title{\normalfont\spacedlowsmallcaps{Levenberg-Marquardt method and partial exact penalty parameter selection}\\ \spacedlowsmallcaps{in bilevel optimization}} 
\author{Andrey Tin$^{\dag}$ and Alain B. Zemkoho$^{\ddag}$\\
{\small{\emph{Centre for Operational Research, Management Sciences and Information Systems (CORMSIS)}}}\\[-1ex]
{\small{\emph{and School of Mathematical Sciences, University of Southampton, SO17 1BJ Southampton, UK}}}}

\date{\today}
\begin{document}

\renewcommand{\sectionmark}[1]{\markright{\spacedlowsmallcaps{#1}}} 
\lehead{\mbox{\llap{\small\thepage\kern1em\color{halfgray} \vline}\color{halfgray}\hspace{0.5em}\rightmark\hfil}} 

\pagestyle{scrheadings}
\maketitle
\setcounter{tocdepth}{2}
\section*{Abstract}
\noindent
We consider the optimistic bilevel optimization problem, known to have a wide range of applications in engineering, that we transform into a single-level optimization problem by means of the lower-level optimal value function reformulation. Subsequently, based on the partial calmness concept, we build an equation system, which is parameterized by the corresponding partial exact penalization parameter. We then design and analyze a Levenberg-Marquardt method to solve this parametric system of equations. Considering the fact that the selection of the partial exact penalization parameter is a critical issue when numerically solving a bilevel optimization problem, we conduct a careful experimental study to this effect, in the context the Levenberg-Marquardt method, while using the Bilevel Optimization LIBrary (BOLIB) series of test problems.

\let\thefootnote\relax\footnotetext{$\dag$ \textit{e-mail: \url{a.tin@soton.ac.uk}.}}

\let\thefootnote\relax\footnotetext{$\ddag$ \textit{e-mail: \url{a.b.zemkoho@soton.ac.uk}.} \\[0.5ex]
AT was supported by a University of Southampton Vice Chancellor Scholarship, while ABZ was
supported by the EPSRC grant EP/V049038/1 and the Alan Turing Institute under
the EPSRC grant EP/N510129/1.}


\section{Introduction}\label{SecIntroduction}

We consider the optimistic bilevel optimization problem
\begin{equation}\label{initialbilev}
\underset{x, \, y}{\min \;}  F(x,y) \; \mbox{ s.t. } \; G(x, y) \leq 0, \;  y\in S(x):= \arg\underset{y}\min~\{f(x,y):\; g(x,y) \leq 0\},\\
\end{equation}
where $F:\mathds{R}^n \times \mathds{R}^m \rightarrow \mathds{R}$, $f:\mathds{R}^n \times \mathds{R}^m \rightarrow \mathds{R}$, $G:\mathds{R}^n \times \mathds{R}^m \rightarrow \mathds{R}^q$, and $g:\mathds{R}^n \times \mathds{R}^m \rightarrow \mathds{R}^p$.
As usual, we refer to $F$ (resp. $f$) as the upper-level (resp. lower-level) objective function and $G$ (resp. $g$) stands for the upper-level (resp. lower-level) constraint function. Solving problem \eqref{initialbilev} is very difficult because of the implicit nature of the lower-level optimal solution set-valued mapping $S : \mathds{R}^n \rightrightarrows \mathds{R}^m$ that partly describes the feasible set of problem \eqref{initialbilev}. Problem   \eqref{initialbilev} has a wide range of applications in engineering; to have a flavour of this, see, e.g., \cite{DempeZemkohoBook} and references therein.

One of the most common single-level transformation of problem \eqref{initialbilev} is the Karush-Kuhn-Tucker (KKT) reformulation, which consists of replacing inclusion $y\in S(x)$ with the equivalent KKT conditions of the lower-level problem, under appropriate assumptions; see, e.g., \cite{Allendestill12, bilevelreform, bilevelmpec10, HerskovitsTanakaLeontiev2013, PinedaByllingMorales2018}. As shown in \cite{bilevelmpec10}, the first drawback of the KKT reformulation is that it is not necessarily equivalent to the original problem \eqref{initialbilev}. Secondly, from the numerical perspective, the KKT reformulation involves derivatives from the lower-level problem that can require the calculation of second (resp. third) order derivatives for first (resp. second) order optimization-type methods. Thirdly, it is shown in the new paper \cite{ZemkohoZhouComparison2021} that the lower-level value function (LLVF) reformulation
\begin{equation}\label{initialvalfuncform0}
\underset{x,\,y}\min~F(x,y) \;\mbox{ s.t. } \; G(x,y)\leq 0, \;\, g(x,y)\leq 0, \;\, f(x,y)\leq \varphi(x),
\end{equation}
where the LLVF $\varphi$ is defined by
\begin{equation}\label{LLVF}
  \varphi(x) := \min~\left\{f(x,y)\left|~g(x,y)\leq 0\right.\right\},
\end{equation}
can lead to a better numerical performance than the KKT reformulation for certain problem classes; in particular for the Bilevel Optimization LIBrary (BOLIB) examples \cite{bolib17}. It is also important note that problem \eqref{initialvalfuncform0}--\eqref{LLVF} is equivalent to the original problem \eqref{initialbilev} without any assumption.

There is a number of recent studies on solution methods for bilevel programs, based on the LLVF reformulation. 
For example, \cite{mitsos08,paulavicius17,wieseman13} develop global optimization techniques for  \eqref{initialvalfuncform0}--\eqref{LLVF}. \cite{lin14,xu14} propose algorithms computing stationary points for \eqref{initialvalfuncform0}--\eqref{LLVF} in the special case where $G$ (resp. $g$) is independent from $y$ (resp. $x$). But, the value function in the latter work is approximated by an entropy function, which is difficult to compute for problems with a big number of lower-level variables, as the corresponding approximation relies on integral calculations.

More recently, \cite{newtonbilevel18, ZemkohoZhouComparison2021} have proposed Newton-type methods for problem \eqref{initialvalfuncform0}--\eqref{LLVF} based on a special transformation that enables the optimality conditions of this problem to be squared. Naturally, detailed optimality conditions for \eqref{initialvalfuncform0}--\eqref{LLVF} are non-square and are directly dealt with in \cite{FTZ20} with Gauss-Newton-type methods. But, these methods require inverse calculations at each iteration for matrices which are not necessarily nonsingular. Hence, to expand the applicability of this method to a wider class of problems, we consider and study the convergence of a regularized version in this paper, commonly known as the Levenberg-Marquardt method.

Recall that the  standard approach to derive optimality conditions for problem \eqref{initialvalfuncform0}--\eqref{LLVF}, necessary to build this Levenberg-Marquardt method, is based on the concept of partial calmness. The consequence of this is that all the methods studied in the papers \cite{FTZ20, newtonbilevel18, ZemkohoZhouComparison2021, lin14,xu14} are based on a partial exact penalization parameter, which needs a special care to be selected. However, none of these papers pays a special attention on how to select this parameter. One of the goals of the current paper is to conduct a numerical study on the best way to select this parameter, based on the BOLIB test set \cite{bolib17}, and in the context of the Levenberg-Marquardt method to be introduced in the next section. This study will certainly enlightening the process of selection the partial exact penalization parameter in a broader class of methods to solve optimistic bilevel optimization problems.

Another difficulty working with the LLVF reformulation is that $\varphi$ is typically non-differentiable function. This will be handled in this paper by using upper estimates of the subdifferential of the function (see, e.g., \cite{dempezemkoho1,newoptcond,bilevelreform,optcondbil95}) that will lead to a relatively simple system of optimality conditions, which depends on the aforementioned partial exact penalization parameter.

The remainder of the paper proceeds as follows. In the next section, we start by transforming the necessary optimality conditions into a system of equations. To do so, we substitute the corresponding complementarity conditions by the well-known Fischer-Burmeister function \cite{fischer92} that we subsequently perturb to build a smooth system of equations.
In Section \ref{PartCalmSec}, we focus our attention on partial exact penalization, including a detailed discussion on how to chose the penalty parameter, as well as the related ill-behaviour.  This analysis plays an important role in the performance of the method, given that a particular attention is paid on ways to avoid the aforementioned ill-behaviour. Furthermore, we will present results on the experimental order of convergence of the method and line search stepsize parameter at the last iteration. The results in Section \ref{SecNumExper} are compared with known solutions of the problems to check the performance of the method under different partial exact penalization parameter scenarios.

\section{The algorithm and convergence analysis}\label{smoothLMsection}
We start this section with some definitions necessary to state optimality conditions for the bilevel program \eqref{initialvalfuncform0}--\eqref{LLVF}, which will be used to construct the algorithm. The lower-level problem is \emph{fully convex} if the functions $f$ and $g_i$, $i=1, \ldots, p$ are convex in $(x,y)$.
A point $(\bar{x},\bar{y})$ is said to be \emph{lower-level regular} if there exists a vector $d\in \mathbb{R}^m$ such that we have
\begin{equation}\label{LMFCQ}
\nabla_y g_i(\bar{x},\bar{y})^\top d < 0 \;\; \mbox{ for } \;\; i\in I_g(\bar{x},\bar{y}):= \left\{i:\;\, g_i (\bar{x}, \bar{y})=0\right\}.
\end{equation}
Obviously, this is the \emph{Mangasarian-Fromovitz constraint qualification} (MFCQ) for the lower-level problem in  \eqref{initialbilev}.
Similarly, a point $(\bar{x},\bar{y})$ is \emph{upper-level regular} if there exists $d\in \mathbb{R}^n\times \mathbb{R}^m$ such that
\begin{equation}\label{UMFCQ}
  \begin{array}{rl}
  \nabla g_i(\bar{x}, \bar{y})^\top d < 0 & \text{for } \;\;  j\in I_g (\bar{x},\bar{y}),\\
  \nabla G_j(\bar{x}, \bar{y})^\top d < 0 & \text{for } \;\; j\in I_G (\bar{x},\bar{y}):=\left\{j:\;\, G_j (\bar{x}, \bar{y})=0\right\}.\\
  \end{array}
\end{equation}
Finally, to write the necessary optimality conditions for problem \eqref{initialvalfuncform0}--\eqref{LLVF}, it is standard to use the following partial calmness concept \cite{optcondbil95}:
\begin{dfn}\label{defpartcalm}
Let $(\bar{x},\bar{y})$ be a local optimal solution of problem \eqref{initialvalfuncform0}--\eqref{LLVF}. The problem is partially calm at $(\bar{x},\bar{y})$ if there exists $\lambda>0$ and a neighbourhood $U$ of $(\bar{x},\bar{y},0)$ such that
$$
F(x,y)-F(\bar{x},\bar{y})+\lambda |u| \geq 0,\;\, \forall (x,y, u)\in U:\; G(x,y)\leq 0, \; g(x,y)\leq 0, \; f(x,y)- \varphi(x)-u=0.
$$
\end{dfn}
The following relationship shows that the partial calmness concept enables the penalization of the optimal value function constraint $f(x,y)-\varphi(x)\leq 0$, to generate a tractable optimization problem, as it is well-known that standard constraint qualifications (including the MFCQ, for example) fail for problem \eqref{initialvalfuncform0}--\eqref{LLVF}; see \cite{dempezemkoho1, optcondbil95}.
\begin{thm}[\cite{optcondbil95}]\label{thresholdlem}
Let $(\bar x, \bar y)$ be a local minimizer of \eqref{initialvalfuncform0}--\eqref{LLVF}. Then, this problem  is partially calm at $(\bar x, \bar y)$ if and only if there is some $\bar \lambda >0$ such that for any $\lambda\geq \bar \lambda$, the point $(\bar x, \bar y)$ is also a local minimizer of
\begin{equation}\label{Penalized Problem}
\underset{x,\,y}\min~F(x,y) + \lambda (f(x,y)-\varphi(x)) \;\mbox{ s.t. } \; G(x,y)\leq 0, \;\,  g(x,y)\leq 0.
\end{equation}
\end{thm}
Problem \eqref{Penalized Problem} is known as the \emph{partial exact penalization} problem of \eqref{initialvalfuncform0}--\eqref{LLVF}, as only the optimal value constraint is penalized. Next, we state the necessary optimality conditions for problem \eqref{initialvalfuncform0}--\eqref{LLVF} based on the aforementioned penalization, while using the upper-and lower-level regularity conditions  and an estimate of the subdifferential of $\varphi$ \eqref{LLVF}; see, e.g., \cite{dempezemkoho1,newoptcond,bilevelreform,optcondbil95}, for the details.

\begin{thm}\label{BilevelOptTh}  Let $(\bar{x},\bar{y})$ be a local optimal solution of problem \eqref{initialvalfuncform0}, where all involved functions are assumed to be continuously differentiable, $\varphi $ is finite around $\bar{x}$, and the lower-level problem is fully convex.
 Furthermore, suppose that the problem is partially calm at $(\bar{x},\bar{y})$, and the lower- and upper-level regularity conditions are both satisfied at $(\bar{x},\bar{y})$. Then, there exist $\lambda > 0$, as well as $u, w\in \mathbb{R}^p$ and $v\in \mathbb{R}^q$ such that
\begin{align}
\nabla F(\bar{x}, \bar{y})+\nabla g(\bar{x}, \bar{y})^T (u - \lambda w)
+ \nabla G(\bar{x},\bar{y})^T v =0, \label{kktbilev12} \\
\nabla_y f(\bar{x}, \bar{y}) + \nabla_y g(\bar{x}, \bar{y})^T w = 0, \label{kktbilev32} \\
u\geq 0, \;\; g(\bar{x}, \bar{y})\leq 0, \;\; u^T g(\bar{x}, \bar{y})=0, \label{kktbilev42} \\
v\geq 0, \;\; G(\bar{x}, \bar{y})\leq 0, \;\; v^T G(\bar{x},\bar{y})=0, \label{kktbilev52} \\
w\geq 0, \;\; g(\bar{x}, \bar{y})\leq 0, \;\; w^T g(\bar{x}, \bar{y})=0. \label{kktbilev62}
\end{align}
\end{thm}
In this result, partial calmness and full convexity are essential and fundamentally related to the nature of the bilevel optimization. Hence, it is important to highlight a few classes of problems satisfying these assumptions. Partial calmness has been the main tool to derive optimality conditions for \eqref{initialbilev} from the perspective of the optimal value function; see, e.g.,  \cite{newoptcond, bilevelreform, dempezemkoho1, optcondbil95}.
It  automatically holds if $G$ is independent from $y$ and the lower-level problem is defined by
\begin{equation}\label{LinearProblems}
f(x,y):= c^\top y \;\; \mbox{ and }\;\, g(x,y):=A(x) + By,
\end{equation}
where $A: \mathbb{R}^n \rightarrow \mathbb{R}^p$,  $c\in \mathbb{R}^m$, and $B\in \mathbb{R}^{p\times m}$.
More generally, various sufficient conditions ensuring that partial calmness holds have been studied in the literature; see \cite{optcondbil95} for the seminal work on the subject. More recently, the paper \cite{MMZ20} has revisited the condition, proposed a fresh perspective, and established new  sufficient conditions for it to be satisfied.

As for full convexity, it will automatically hold for the class of problem defined in \eqref{LinearProblems}, provided that each component of the function $A$ is convex. Another class of nonlinear fully convex lower-level problem is given in \cite{LamparielloSagratellaNumerically2020}. Note however that when it is not possible to guarantee that this assumption is satisfied, there are at least two alternative scenarios to obtain the same optimality conditions as in Theorem \ref{BilevelOptTh}. The first is to replace the full convexity assumption by the \emph{inner semicontinuity} of the optimal solution set-valued mapping $S$ \eqref{initialbilev}. Secondly, note that a much weaker qualification condition known as \emph{inner semicompactness} can also be used here. However, under the latter assumption, it will additionally be required to have $S(\bar x)=\{\bar y\}$ in order to get the optimality conditions in \eqref{kktbilev12}--\eqref{kktbilev62}. The concept of inner semicontinuity (resp. semicompactness) of $S$ is closely related to the lower semicontinuity (resp. upper semicontinuity) of set-valued mappings; for more details on these notions and their ramifications on bilevel programs, see, e.g., \cite{newoptcond,bilevelreform,dempezemkoho1}. 

It is important to mention that various other necessary optimality conditions, different from the above ones, can be obtained, depending on the assumptions made. Details of different stationarity concepts for \eqref{initialvalfuncform0} can be found in the latter references, as well as in \cite{zemkohothesis}.

To reformulate the complementarity conditions \eqref{kktbilev42}--\eqref{kktbilev62} into a system of equations, we use the well-known  \emph{Fischer-Burmeister} function \cite{fischer92} $\phi(a,b) := \sqrt{a^2+b^2}-a-b$, which ensures that
\[
\phi(a,b)=0\;\; \iff \;\; a\geq 0, \;\; b \geq 0, \;\; ab=0.
\]
This leads to the reformulation of the optimality conditions \eqref{kktbilev12}--\eqref{kktbilev62} into the system of equations:
\begin{align}
\Upsilon^\lambda (z) := \left(\begin{array}{rr}  \nabla_x F(x, y)+\nabla_x g(x, y)^T (u - \lambda  w)
+ \nabla_x G(x,y)^T v \\
  \nabla_y F(x, y) +\nabla_y g(x, y)^T (u-\lambda w)+ \nabla_y  G(x,y)^T v\\
  \nabla_y f(x, y) + \nabla_y g(x, y)^T w \\
 \sqrt{u^2+g(x, y)^2} - u+g(x, y) \\
  \sqrt{v^2+G(x,y)^2 } - v+G(x,y) \\
  \sqrt{w^2+g(x, y)^2 } - w+g(x, y) \end{array} \right) =
0, \label{bilevncp21}
\end{align}
where we have $z:=(x, y, u, v, w)$ and
\begin{equation}\label{hunam}
\sqrt{u^2+g(x, y)^2} - u+g(x, y) :=
\left(\begin{array}{c}
\sqrt{u_1^2+g_1(x, y)^2} - u_1+g_1(x, y) \\
\vdots \\
\sqrt{u_p^2+g_p(x, y)^2} - u_p+g_p(x, y)
 \end{array} \right).
\end{equation}
$\sqrt{v^2+G(x,y)^2 } - v+G(x,y)$ and $ \sqrt{w^2+g(x, y)^2 } - w+g(x, y)$ are defined as in \eqref{hunam}. The superscript $\lambda$ is used to emphasize the fact that this number is a parameter and not a variable for equation \eqref{bilevncp21}. One can easily check that this system is made of $n+2m+p+q+ p$ real-valued equations and $n+m+p+q+p$ variables. Clearly, this means that \eqref{bilevncp21} is an \emph{overdetermined} system and the Jacobian of $\Upsilon^\lambda (z)$, when it exists, is a non-square matrix.

To focus our attention on the main ideas of this paper, we smoothen the function $\Upsilon^\lambda$ \eqref{bilevncp21} using the \emph{smoothed  Fischer-Burmeister function} (see \cite{kanzow1996}) defined by
\begin{equation} \label{FischerBurmsmooth}
\phi^\mu_j (x, y, u) :=\sqrt{u^2_j+g_j(x, y)^2 + 2\mu} - u_j+g_j(x, y), \; j=1, \ldots, p,
\end{equation}
where the perturbation $\mu >0$ helps to guaranty its differentiability at points $(x, y, u)$, where $u_j = g_j(x,y)=0$  for $j=1, \ldots, p$. It is well-known (see latter reference) that for $j=1, \ldots, p$,
\begin{equation}\label{equivalent-ce}
\phi^\mu_j (x, y, u) = 0 \;\; \Longleftrightarrow \;\; \left[ u_j>0,\; -g_j (x,y) >0, \; -u_j g_j(x,y) = \mu\right].
\end{equation}
The smoothed version of system \eqref{bilevncp21} then becomes
\begin{align}
\Upsilon^{\lambda}_\mu (z)=  \left(\begin{array}{rr}  \nabla_x F(x, y)+\nabla_x g(x, y)^T (u - \lambda  w)
+ \nabla_x G(x,y)^T v \\
  \nabla_y F(x, y) +\nabla_y g(x, y)^T (u-\lambda w)+ \nabla_y  G(x,y)^T v\\
  \nabla_y f(x, y) + \nabla_y g(x, y)^T w \\
 \sqrt{u^2+g(x, y)^2 + 2\mu} - u+g(x, y) \\
  \sqrt{v^2+G(x,y)^2 + 2\mu } - v+G(x,y) \\
  \sqrt{w^2+g(x, y)^2  + 2\mu} - w+g(x, y) \end{array} \right) =
0, \label{bilevncp2222}
\end{align}
following the convention in \eqref{hunam}, where $\mu$ is a vector of appropriate dimensions with sufficiently small positive elements.
Under the assumption that all the functions involved in problem \eqref{initialbilev} are continuously differentiable, $\Upsilon^{\lambda}_\mu$ is also a continuously differentiable function for any $\lambda > 0$ and $\mu>0$. We can easily check that for a fixed value of $\lambda > 0$,
\begin{equation}\label{SmoothingScheme}
\|\Upsilon^{\lambda}_{\mu} (z) - \Upsilon^{\lambda} (z)\| \;\; \longrightarrow 0 \;\; \mbox{ as } \;\; \mu \downarrow 0.
\end{equation}
Hence, based on this scheme,
our aim is to consider a sequence $\{\mu_k\}$ decreasing to $0$ such that equation \eqref{bilevncp21} is approximately solved while leading to
\[
\Upsilon^{\lambda}_{\mu^k} (z^k)  \;\; \longrightarrow 0 \;\; \mbox{ as } \;\; k \uparrow \infty
\]
for a fixed value of $\lambda >0$. In order to proceed, let us define the least squares problem
\begin{equation}\label{Lsquares}
\min~\Phi^\lambda_{\mu}(z) := \frac{1}{2} \norm{\Upsilon^\lambda(z)}^2.
\end{equation}
Before we introduce the smoothed Levenberg-Marquardt method that will be one of the main focus points of this paper, note that for fixed $\lambda >0$ and $\mu>0$,  
\begin{equation} \label{compactJsmooth}
\nabla \Upsilon^\lambda_\mu (z) = \left[
  \begin{array}{cccc}
    \nabla^2 L^\lambda (z)  & \nabla g(x,y)^T & \nabla G(x,y)^T & -\lambda \nabla g(x,y)^T \\
     \nabla (\nabla_y\mathcal{L} (z)) & O & O & \nabla_y g(x,y)^T \\
    \mathcal{T}^\mu \nabla g(x,y) & \Gamma^\mu & O & O \\
     \mathcal{A}^\mu \nabla G(x,y) & O &  \mathcal{B}^\mu & O \\
    \Theta^\mu \nabla g(x,y) & O & O & \mathcal{K}^\mu \\
  \end{array}
\right]
\end{equation}
with  
the pair $\left(\mathcal{T}^\mu, \Gamma^\mu\right)$ defined by $\mathcal{T}^\mu := diag~\{\tau_1^\mu,..,\tau_p^\mu\}$ and $\Gamma^\mu := diag~\{\gamma_1^\mu,..,\gamma_p^\mu\}$, where
\begin{equation}\label{deftau-mu}
\tau^{\mu}_j:= \frac{g_j (x,y)}{\sqrt{u_j^2 + g_j(x,y)^2 + 2\mu}}+1 \; \mbox{ and } \; \gamma^{\mu}_j:= \frac{u_j}{\sqrt{u_j^2 + g_j(x,y)^2 +2\mu}} - 1, \;\, j = 1, \ldots p.
\end{equation}
The pairs  $\left(\mathcal{A}^{\mu},\, \mathcal{B}^{\mu}\right)$ and $\left(\Theta^{\mu}, \mathcal{K}^{\mu}\right)$ are defined in a similar way, based on $(v_j , G_j(x,y))$, $j=1,\ldots,q$ and $(w_j,g_j(x,y))$, $j=1,\ldots,p$, respectively.

We now move on to present some definitions that will help us state the algorithm with \emph{line search}. It is well-known that line search helps to choose the optimal step length to avoid over-going an optimal solution in the direction $d^k$ and also to globalize the convergence of the method, i.e., have more flexibility on the starting point $z^0$.
The optimal step length $\gamma_k$ can be calculated through minimizing $\Phi^\lambda(z^k + \gamma_k d^k)$, with respect to $\gamma_k$, such that
\[
\Phi^\lambda(z^k + \gamma_k d^k) \leq \Phi^\lambda(z^k) + \sigma \gamma_k \nabla \Phi^\lambda (z^k)^T d^k \; \mbox{ for } \; 0<\sigma<1.
\]
That is, we are looking for $\gamma_k = arg \min_{\gamma \in \mathbb{R}} ||\Upsilon^\lambda (z^k + \gamma d^k)||^2$.
To implement line search, it is standard to use \emph{Armijo condition} that guarantees a decrease at the next iterate.
\begin{dfn}\label{Line search 1}Fixing $d$ and $z$, consider the function $\phi_\lambda (\gamma) := \Phi^\lambda(z+\gamma d)$.  Then,  the \emph{Armijo} condition  will be said to hold if $\phi_\lambda(\gamma)\leq \phi(0) + \gamma \sigma \phi'_\lambda(\gamma)$ for some $0<\sigma <1$.
\end{dfn}
The practical implementation of the Armijo condition is based on backtracking.
\begin{dfn}\label{Line search 2}
Let $\rho \in (0,1)$ and $\bar{\gamma} >0$. \emph{Backtracking} is the process of checking over a sequence $\bar{\gamma}$, $\rho \bar{\gamma}$,  $\rho^2 \bar{\gamma}$, \ldots, until a number $\gamma$ is found satisfying the Armijo condition.
\end{dfn}

Line search is widely used in continuous optimization; see, e.g., \cite{nocedal99} for more details. For the implementation in this paper, we start with stepsize $\gamma_0:=1$; then, if the condition
\[
\norm{\Upsilon^{\lambda}(z^k+\gamma_k d^k)}^2 < \norm{\Upsilon^{\lambda}(z^k)}^2 + \sigma \gamma_k \nabla \Upsilon_\mu^\lambda(z^k)^T \Upsilon^\lambda(z^k) d^k,
\]
is not satisfied, we set $\gamma_k = \gamma_k/2$ and check again until the condition above is satisfied.
More generally, the algorithm proceeds as follows:
 \begin{alg}
Smoothed Levenberg-Marquardt Method for Bilevel Optimization
\label{algorithm LMLs}
\begin{algorithmic}
 \STATE \textbf{Step 0}: Choose $(\lambda, \mu, K, \epsilon, \alpha_0)\in \left(\mathbb{R}^*_+\right)^5$, $(\rho, \sigma, \gamma_0) \in (0,1)^3$,  $z^0:=(x^0, y^0, u^0, v^0, w^0)$,  and set $k:=0$.  
 \STATE \textbf{Step 1}: If $\norm{\Upsilon^{\lambda}_\mu (z^k)}<\epsilon$  or $k\geq K$, then stop.
 \STATE \textbf{Step 2}: Calculate the Jacobian $\nabla \Upsilon_\mu^\lambda(z^k)$ and subsequently find a vector $d^k$ satisfying 
 \begin{equation}\label{SLMdirnum}
\left(\nabla \Upsilon_\mu^\lambda (z^k)^\top \nabla  \Upsilon_\mu^\lambda(z^k) + \alpha_k I\right) d^k = - \nabla  \Upsilon_\mu^\lambda(z^k)^\top  \Upsilon^\lambda_\mu(z^k),
\end{equation}
where $I$ denotes the identity matrix of appropriate size.
\\
 \STATE \textbf{Step 3}: \textbf{While} $\norm{\Upsilon^\lambda_\mu(z^k+\gamma_k d^k)}^2 \geq \norm{\Upsilon^{\lambda}_\mu(z^k)}^2 + \sigma \gamma_k \nabla \Upsilon_\mu^\lambda(z^k)^T \Upsilon^\lambda_\mu(z^k) d^k$, \textbf{do} $\gamma_k \leftarrow \rho \gamma_k$ \textbf{end}.
\STATE \textbf{Step 4}: Set $z^{k+1}:=z^k + \gamma_k d^k$,  $k:=k+1$, and go to Step 1.
\end{algorithmic}
\end{alg}

Note that in Step 0, $\mathbb{R}^*_+:=(0, \, \infty)$.
Before we move on to focus our attention on the practical implementation details of this algorithm, we present its convergence result, which is based on the following selection of the Levenberg–Marquardt (LM) parameter $\alpha_k$:
\begin{equation}\label{LM parameter}
\alpha_k = \|\Upsilon_\mu^\lambda(z^k)\|^{\eta} \,\mbox{ for any choice of }\, \eta\in \left[1, \, 2\right].
\end{equation}
\begin{thm}[\cite{FanYuan05}]  \label{LMconvtheoremEBs1}
Consider Algorithm \ref{algorithm LMLs} with fixed values for the parameters $\lambda>0$ and $\mu>0$ and let $\alpha_k$ be defined as in \eqref{LM parameter}. 
Then, the sequence $\{z^k\}$ generated by the algorithm converges quadratically to $\bar{z}$ satisfying $\Upsilon_\mu^\lambda(\bar z)=0$, under the following assumptions:
\begin{enumerate}
\item  $\Upsilon_\mu^\lambda : \mathbb{R}^{N} \rightarrow \mathbb{R}^{N+m}$ is continuously differentiable and $\nabla \Upsilon_\mu^\lambda : \mathbb{R}^{N} \rightarrow \mathbb{R}^{(N+m)\times(N)}$ is locally Lipschitz continuous in a neighbourhood of $\bar z$.
\item There exists some $C>0$ and $\delta>0$ such that
\[
C dist(z,Z^\lambda_\mu) \leq \norm{\Upsilon_\mu^\lambda (z)} \phantom{-} \text{for all }\;\, z\in \mathcal{B} (\bar{z}, \delta),
\]
where \emph{dist} denotes the distance function and $Z^\lambda_\mu$ corresponds to the solution set of equation \eqref{bilevncp2222}.
\end{enumerate}
\end{thm}
For fixed values of $\lambda>0$ and $\mu>0$, assumption 1 in this theorem is automatically satisfied if all the functions involved in problem \eqref{initialbilev} are twice continuously differentiable. According to \cite{YF01},  assumption 2 of Theorem \ref{LMconvtheoremEBs1} is fulfilled if the matrix $\nabla \Upsilon_\mu^\lambda$ has a full column rank. Various conditions guarantying that $\nabla \Upsilon_\mu^\lambda$ has a full rank have been developed in \cite{FTZ20}.
Below, we present an example of bilevel program satisfying the first and second assumptions of Theorem \ref{LMconvtheoremEBs1}. 
\begin{exm}\label{ExampIndepAssumSmooth}
Consider the following instance of problem \eqref{initialbilev} from the BOLIB library \cite[LamprielloSagratelli2017Ex33]{bolib17}:
\[
\begin{array}{rll}
  F(x,y)& := & x^2 + (y_1+y_2)^2,\\	
  G(x,y)& := & -x+0.5,\\
  f(x,y)& := &y_1,\\
  g(x,y)& := & \left(\begin{array}{c}  -x-y_1-y_2+1\\
-y  \end{array} \right).
\end{array}
\]
Obviously, assumption 1 holds.
According to \cite{FTZ20}, for this example, the columns of $\nabla \Upsilon_\mu^\lambda$ are linearly independent at the solution point
\[
\bar z := (\bar{x}, \,\bar{y}_1, \, \bar{y}_2, \, \bar{u}_1, \, \bar{u}_2, \, \bar{u}_3, \, \bar{v}, \, \bar{w}_1,\, \bar{w}_2, \, \bar{w}_3)=(0.5, \, 0, \, 0.5, \, 1, \, \lambda, \, 0, \, 0, \, 0, \, 1, \, 0)
\]
with the parameters chosen as $\mu=2\times10^{-2}$ and $\lambda=10^{-2}$.
\qed
\end{exm}

On the selection of the LM parameter $\alpha_k$, we conducted a preliminary analysis based on
the BOLIB library test set \cite{bolib17}. It was observed that for almost all the corresponding examples, the choice $\alpha_k := \left\|\Upsilon_\mu^\lambda(z^k)\right\|^\eta$ with $\eta \in (1, 2]$ leads to a very poor performance of  Algorithm \ref{algorithm LMLs}. The typical behaviour of the algorithm for  $\eta \in (1, 2]$  is shown in the following example.

\begin{exm}Consider the following scenario of problem \eqref{initialbilev} from \cite[AllendeStill2013]{bolib17}:
\[
\begin{array}{rll}
  F(x,y)& := & x_1^2 - 2x_1 + x_2^2 - 2x_2 + y_1^2 + y_2^2,\\	
  G(x,y)& := & \left(\begin{array}{c}
                        -x\\
                        -y\\
                        x_1-2
                     \end{array}
     \right),\\
  f(x,y)& := &y_1^2 - 2x_{1}y_1 + y_2^2 - 2x_{2}y_2,\\
  g(x,y)& := & \left(\begin{array}{c}
                        (y_1 - 1)^2 -0.25 \\
                        (y_2 - 1)^2 -0.25
                     \end{array}
     \right).
\end{array}
\]
Figure \ref{alpha2} shows the progression of $\norm{\Upsilon^{\lambda}(z^k)}$ generated from Algorithm \ref{algorithm LMLs} with $\alpha_k$ selected as in \eqref{LM parameter} while setting $\eta=1$ and $\eta=2$, respectively. Clearly, after about 100 iterations $\norm{\Upsilon^{\lambda}(z^k)}$ blows up relatively quickly when $\eta=2$, while it falls and stabilizes within a certain tolerance for $\eta=1$.
\begin{figure}[H]
    \centering
   \subfloat[$\alpha_k:=\norm{\Upsilon^{\lambda}(z^k)}$]{{\includegraphics[width=7.25cm]{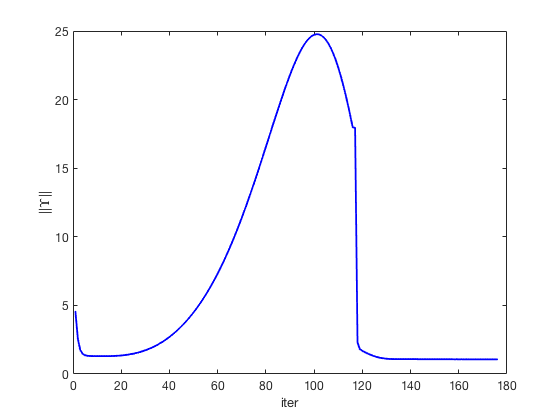} }}%
           \qquad
   \subfloat[$\alpha_k:=\norm{\Upsilon^{\lambda}(z^k)}^2$]{{\includegraphics[width=7.25cm]{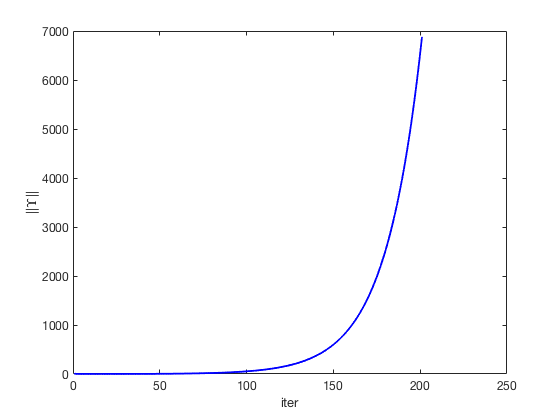} }}%
    \caption{Typical behaviour of Algorithm \ref{algorithm LMLs} for two scenarios of the LM parameter}%
    \label{alpha2}
\end{figure}
It is worth noting that the scale on the y-axis of Figure \ref{alpha2}(b) is quite large. Hence, it might not be apparent that solutions at the early iterations of the algorithm are much better for $\eta=1$ compared to the ones obtained in the case where $\eta=2$.  \hfill \qed
\end{exm}
Note that for almost all of the examples in the BOLIB test set \cite{bolib17}, we have a behaviour of Algorithm \ref{algorithm LMLs} similar to that of Figure \ref{alpha2}(b) when $\alpha_k=\norm{\Upsilon^{\lambda}(z^k)}^\eta$ for many different choices of $\eta\in(1,2]$.
It is important to mention that the behaviour of the algorithm is not surprising, as it is well-known in the literature that with the sequence
$\alpha_k:=\norm{\Upsilon_\mu^\lambda(z^k)}^2$, for example, the Levenberg–Marquardt method often faces some issues. Namely, when the sequence $z^k$ is close to the solution set, $\alpha_k:=\norm{\Upsilon_\mu^\lambda(z^k)}^2$  could become smaller than the machine precision and hence lose its role as a result of this. On the other hand, when the sequence is far away from the solution set, $\norm{\Upsilon^\lambda_\mu(z^k)}^2$  may be very large; making a movement towards the solution set to be very slow. Hence, from now on, we use  $\alpha_k:=\norm{\Upsilon^{\lambda}(z^k)}$,  with $k=0, \, 1, \ldots$, for all the analysis of Algorithm \ref{algorithm LMLs} conducted in this paper. Note however that there are various other approaches to select the LM parameter $\alpha_k$ in the literature; see, e.g., \cite{behling16, FanYuan05, YF01}.

To conclude this section, it is important to recall that from the perspective of bilevel optimization, the partial exact penalization parameter $\lambda$ is a key element of Algorithm \ref{algorithm LMLs}, as it originates from the penalization of the value function constraint $f(x,y)\leq \varphi(x)$. Additionally, unlike the other parameters involved in the algorithm, which have benefited from many years of research, as it will be clear in Section \ref{SecNumExper}, it remains unknown what is the best way to select $\lambda$ while solving problem \eqref{initialbilev} via the value function reformulation \eqref{initialvalfuncform0}--\eqref{LLVF}. Hence, the focus of the remaining parts of this paper will be on the selection of $\lambda$ and assessing its  impact on the efficiency of Algorithm \ref{algorithm LMLs}.

\section{Partial exact penalty parameter selection}\label{PartCalmSec}
The aim of this section is to explore the best way to select the penalization parameter $\lambda$. Based on Theorem \ref{thresholdlem}, we should be getting the solution for some threshold value $\bar \lambda$ of the penalty parameter and the algorithm should be returning the value of the solution for any $\lambda \geq \bar \lambda$. Hence, increasing $\lambda$ at every iteration seems to be the ideal approach to follow this logic to obtain and retain the  solution. Hence, to proceed, we use the increasing sequence $\lambda_k := 0.5*(1.05)^k$, where $k$ is the number of iterations of the algorithm. The main reason of this choice is that the final value of $\lambda$ need not to be too small to recover solutions and not too large to avoid ill-conditioning. It was observed that going too aggressive with the sequence (e.g., with $\lambda_k:=2^k$) forces the algorithm to diverge.
Also, it was observed that for fixed and small values of $\lambda$ (i.e., $\lambda <1$), Algorithm \ref{algorithm LMLs} performs well for many examples.  This justifies choosing the starting value for varying parameter less than $1$. Overall, the aim here is vary $\lambda$ as mentioned above and assess what could be the potential best ranges for selection of the parameter based on our test set from  BOLIB \cite{bolib17}.

\subsection{How far can we go with the value of $\lambda$?}\label{illbehaveSec}
We start here by acknowledging that it is very difficult to check that partial calmness (cf. Definition \ref{defpartcalm}) holds in practice. Nevertheless, we would like to consider Theorem \ref{thresholdlem} as the base of our experimental analysis, and ask ourselves how large does $\lambda$ need to be for Algorithm \ref{algorithm LMLs} to converge. Intuitively, one would think that taking $\lambda$ as whatever large number should be fine.
However, this is usually not the case in practice. One of the main reasons for this is that for too large values of $\lambda$, Algorithm \ref{algorithm LMLs} does not behave well. In particular, if we run Algorithm \ref{algorithm LMLs}  with varying $\lambda$ for a very large number of iterations, the  value of the $Error$ blows up  at some point and the values $\norm{\Upsilon^{\lambda}_\mu (z^k)}$ stop decreasing.
Recall that ill-conditioning (and possibly the \emph{ill-behaviour}) refers to the fact that one eigenvalue of the Hessian involved in $\nabla \Upsilon_\mu^\lambda(z^k)$ being much larger than the other eigenvalues, which affects the curvature in the negative way for gradient methods \cite{PTVF92}.   
To analyze this ill-behaviour in this section, we are going to run our algorithm for $1,000$ iterations with no stopping criteria and let $\lambda$ vary indefinitely. We will subsequently look at which iteration the algorithm blows up and record the value of $\lambda$ there. To proceed, we denote by $\lambda_{ill}$ the first value of $\lambda$ for which the ill-behaviour is observed for each example. The table below presents the number of BOLIB examples that are approximately within certain intervals of $\lambda_{ill}$.
\begin{small}
\begin{longtable}{c|c|c|c|c|c}
\caption{Approximate intervals for the values of $\lambda$ when the ill-behaviour starts}\\[2ex]
  \hline
$\lambda_{ill}$ & $\lambda_{ill}<10^7$ & $10^7<\lambda_{ill}<10^9$ & $10^9<\lambda_{ill}<10^{11}$ & $10^{11}<\lambda_{ill}<10^{20}$ & Not observed  \\
\hline
Examples&1&6&72&11&34\\\hline
\end{longtable}
\end{small}
On average, the ill-behaviour seems to typically occur after about $500$ iterations (where $\lambda_{ill}\approx10^{10}$), as it can be seen in the table above. For $34$ problems ill behaviour was not observed under the scope of $1000$ iterations. We also see that for most of the problems ($72/124$), the ill-behaviour happens for $10^9<\lambda<10^{11}$.
We further observe that
algorithm has shown to behave well for the values of penalty parameter  $\lambda<10^9$ with only $7/124$ examples demonstrating ill behaviour for such $\lambda$.
This makes the choice of very large values of $\lambda$ not attractive at all.
Mainly, the analysis shows that choosing $\lambda>10^9$ could cause the algorithm to diverge. Hence, based on our method, selecting $\lambda\leq10^7$ seems to be a very safe choice for our BOLIB test set. This is useful for the choice of fixed $\lambda$ as we can choose values smaller than $10^7$.
For the case of varying $\lambda$ the values are controlled by the stopping criteria proposed in the next section.
The complete results on the values of $\lambda_{ill}$ for each example will be presented in Table \ref{combinedtable} below.


Interestingly,  34 out of the 124 test problems do not demonstrate any ill behaviour signs even if we run the algorithm for $1,000$ iterations with $\lambda=0.5*1.05^{iter}$. A potential reason for this could just be that the parameter $\lambda$ does not get  sufficiently large after $1,000$ iterations to cause problems for these problems. It could also
be that the eigenvalues of the Hessian are not affected by large values of $\lambda$ for these examples. Also, there is a possibility that  elements of the Hessian involved in \eqref{SLMdirnum} do not depend on $\lambda$ at all, as for $20/34$ problems where the function $g$ is linear in $(x,y)$ or not present in these problems.
Next, we focus on assessing what magnitudes of the penalty parameter $\lambda$ seem to lead to the bet performance for our method.

\subsection{Do the values of $\lambda$ really need to be large?}\label{Seclamvals}
It is clear from the previous subsection that to reduce the chances for Algorithm \ref{algorithm LMLs} to diverge or exhibit some ill-behaviour, we should approximately select $\lambda < 10^7$. However, it is still unclear whether only large values of $\lambda$ would be sensible to ensure a good behaviour of the algorithm. In other words, it is important to know whether relatively small values of $\lambda$ can lead to a good behaviour for Algorithm \ref{algorithm LMLs}. To assess this, we attempt here to identify inflection points, i.e., values of $k$ where we have $\norm{\Upsilon^{\lambda_k}_\mu (z^k)}<\epsilon$ as $\lambda_k$ varies increasingly as described in the introduction of this section.
We would then record the value of $\lambda_k$ at these points.
%
%

\begin{figure}[H]
    \centering
    \subfloat[Complete run 1]{{\includegraphics[width=4.2cm]{LMdraw2.png} }}%
   \subfloat[Inflection point stop 1]{{\includegraphics[width=4.2cm]{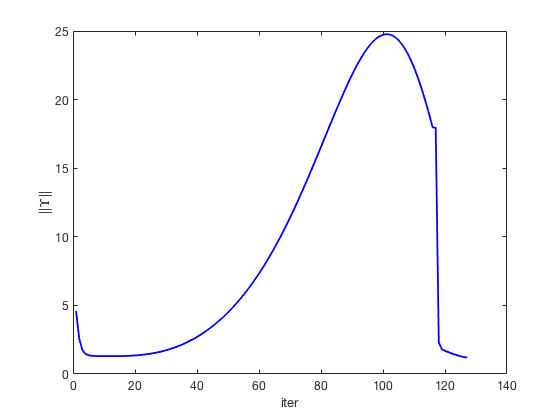} }}%
   \subfloat[Complete run 2]{{\includegraphics[width=4.2cm]{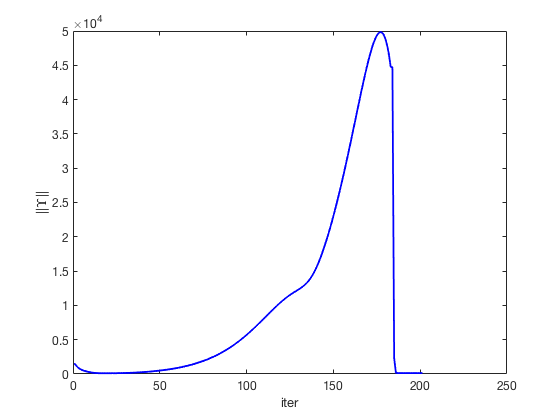}}}%
   \subfloat[Inflection point stop 2]{{\includegraphics[width=4.2cm]{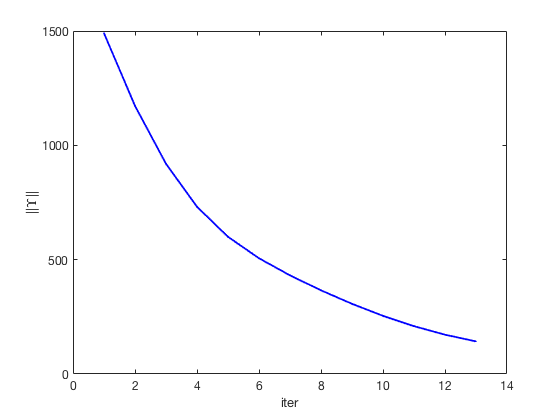}}}%
    \caption{(a) and (b) illustrate the inflection point identification for Example \texttt{AllendeStill2013} and similarly (c) and (d) correspond to Example  \texttt{Anetal2009}; the examples are taken from BOLIB \cite{bolib17}.}%
    \label{twofigsexample2}%
\end{figure}

Ideally, we want to get the threshold $\bar{\lambda}$ such that solution is retained for all $\lambda>\bar \lambda$ in the sense of Theorem \ref{thresholdlem}. To proceed, we extract the information on the final $Error^*:=\norm{\Upsilon^\lambda(z^*)}$ for each example from \cite{LMnumer20} and then rerun the algorithm with varying penalty parameter $\lambda_k:=0.5*1.05^k$ with new stopping criterion (i.e., $Error \leq 1.1 Error^*$) while relaxing all of the stopping criteria (see details in Section \ref{SecNumExper}). This way, we stop once we observe $Error_k:=\norm{\Upsilon^{\lambda_k}(z)}$ close to the $Error^*$ that we obtained in our experiments  \cite{LMnumer20}.
It is worth noting that it would make sense to test only 72/117 ($61.54\%$) of examples, for which algorithm performed well and produced a good solution.
This approach can be thought of as finding the inflection point. For instance, if we have an algorithm running as below, we want to stop at the inflection point after 125-130 iterations. The illustration of how we aim to stop at the inflection point is presented in  Figure \ref{twofigsexample2} (a)--(b) below, where we have $\norm{\Upsilon^\lambda(z)}$ on the y-axis and iterations on the x-axis.

For some of the examples, we obtained a better $Error$ than initial $Error^*$. For these cases, we stopped very early as $Error_k \leq1.1Error^*$ was typically met at an early iteration $k$ (where $\lambda$ was still small), as demonstrated in Figure \ref{twofigsexample2}(c)--(d) above.
This demonstrates the disadvantage of $\lambda$ being an increasing sequence. If the algorithm makes many iterations, the parameter $\lambda$ keeps increasing without possibility to go back to the smaller values. It turns out that  for some examples the smaller value of $\lambda$ was as good as large values, or even better to recover a solution. 
This further justifies the choice to start from the small value of $\lambda$, that is $\lambda_0<1$ and increase it slowly.

With the setting to stop whenever $Error_k \leq1.1 Error^*$,  we often stop very early. Hence, we do not always get $\bar \lambda$ that represents the inflection point, which we aimed to get; cf. Figure \ref{twofigsexample2}(c)--(d), where (c)  has a scale of $10^4$ on the y-axis.
Although, we can clearly see from Figure \ref{twofigsexample2} (c) that inflection point lies around 190-200 iterations, where the value of $\lambda$ is much bigger. It is clear that we stopped earlier due to having small value of $Error$ after 12-45 iterations.
It was observed that such scenario is typical for the examples in the considered test set.  
For this reason, we want to introduce $\lambda^*$ as the \emph{large threshold} of $\lambda$. This value will be used to represent the value of the penalty parameter at the inflection point, where solution starts to be recovered for large values of $\lambda$ ($\lambda>6.02$),  while we also record $\bar \lambda$ as the first (smallest) value of $\lambda$ for which good solution was obtained. 
For instance, with $\lambda$ defined as  $\lambda := 0.5\times1.05^k$ in Figure \ref{twofigsexample2} (c) we have $\bar \lambda = 0.5\times1.05^{12}$ and $\lambda^*=0.5\times1.05^{190}$ as we obtain good solution for small $\lambda$ after 12 iterations and for large $\lambda$ after 190 iterations.
We shall note that the value $\lambda>6.02$ is the value of the penalty parameter after we make at least 50 iterations as for the case with varying $\lambda$ we have $\lambda= 0.5\times1.05^{51}=6.02$.

The complete results of detecting $\bar \lambda$ and $\lambda^*$ is presented in Table \ref{combinedtable} below.
It was observed that the behaviour of the method follows the same pattern for the majority of the examples.
Typically, we get a good solution retaining for a few small values of $\lambda$, then value of the Error blows up and takes some iterations to start decreasing, coming back to obtaining and retaining a good solution for large values of $\lambda$.
Such pattern is clearly demonstrated in Figure \ref{twofigsexample2} (c).
Such a behaviour is interesting as usually, only large values of penalty parameters are used in practice \cite{burke91,pillo89}, as intuitively suggested by Theorem \ref{thresholdlem}. As mentioned in \cite{Fletcher75}, some methods require penalty parameters to increase to infinity to obtain convergence.

\begin{footnotesize}
\begin{longtable}{r|l|c|cHH|cHH|c|l}
\caption{Ill behaviour and two thresholds for $\lambda$ }\\
  \hline
  \textbf{Problem number} &   \textbf{Problem name} &    $\textbf{iter}_{ill} $& $\lambda_{ill}$ &  \textbf{Problem number} & \textbf{Problem name} & $\bar \lambda$ & $Error_{new}$ & $Error_{old}$ & $\lambda^*$ &   iter for $\lambda^*$
 \\
\hline
\rowcolor[gray]{.9}
1 & \texttt{AiyoshiShimizu1984Ex2}&495&1.54e+10 & 1 & \texttt{AiyoshiShimizu1984Ex2}&$9.92*10^4$&1.78&1.65 &$9.92*10^4$&250\\\hline
2 & \texttt{AllendeStill2013} &Not observed&NA & 2 & \texttt{AllendeStill2013} &245&1.24&1.07&245&127\\\rowcolor[gray]{.9}
3 &\texttt{AnEtal2009}&Not observed&NA & - & - & - & - & - & - & -\\\hline
4 & \texttt{Bard1988Ex1}	&520&5.22e+10 & 4 & \texttt{Bard1988Ex1}&0.608&3.53&3.11&245&127\\\rowcolor[gray]{.9}
5 & \texttt{Bard1988Ex2}&536&1.14e+11 & - & - & - & - & - & - & - \\\hline
6 & \texttt{Bard1988Ex3}&Not observed&NA & 6 & \texttt{Bard1988Ex3}&0.525&6.16&7.91&54.1&96\\\rowcolor[gray]{.9}
7 & \texttt{Bard1991Ex1}	&Not observed&NA & 7 & \texttt{Bard1991Ex1}&0.608&2.9&2.52&183&121\\\hline
8 & \texttt{BardBook1998}&502&2.17e+10 & - & - & - & - & - & - & - \\\rowcolor[gray]{.9}
9 & \texttt{CalamaiVicente1994a}&476&6.1e+09 & 9 & \texttt{CalamaiVicente1994a}&0.608&0.573&0.487&27.3&82\\\hline
10 & \texttt{CalamaiVicente1994b}&490&1.21e+10 & 10 & \texttt{CalamaiVicente1994b}&0.739&0.557&0.489&79.9&104\\\rowcolor[gray]{.9}
11 & \texttt{CalamaiVicente1994c}&490&1.21e+10 & - & - & - & - & - & - & - \\\hline
12 & \texttt{CalveteGale1999P1}&538&1.26e+11 & 12 & \texttt{CalveteGale1999P1}&0.525&35.2&32.6&$1.57*10^3$&165\\\rowcolor[gray]{.9}
13 & \texttt{ClarkWesterberg1990a}&520&5.22e+10 & - & - & - & - & - & - & - \\\hline
14 & \texttt{Colson2002BIPA1}&510&3.2e+10 & 14 & \texttt{Colson2002BIPA1}&792&9.74&8.16&792&151\\\rowcolor[gray]{.9}
15 &  \texttt{Colson2002BIPA2}&900&5.88e+18 & 15 &  \texttt{Colson2002BIPA2}&1.33&2.63&2.2&$1.65*10^3$&166\\\hline
16 & \texttt{Colson2002BIPA3}&110&107 & - & - & - & - & - & - & - \\\rowcolor[gray]{.9}
17 & \texttt{Colson2002BIPA4}	&Not observed&NA & 17 & \texttt{Colson2002BIPA4}&0.551&4.01&5.62&$2.68*10^3$&176\\\hline
18 & \texttt{Colson2002BIPA5}	&550&2.25e+11 & - & - & - & - & - & - & - \\\rowcolor[gray]{.9}
19 & \texttt{Dempe1992a}&Not observed&NA & - & - & - & - & - & - & - \\\hline
20 & \texttt{Dempe1992b}	&Not observed&NA & 20 & \texttt{Dempe1992b}	&0.551&5.61&6.33&$1.29*10^3$&161\\\rowcolor[gray]{.9}
21 & \texttt{DempeDutta2012Ex24}&Not observed&NA & - & - & - & - & - & - & - \\\hline
22 & \texttt{DempeDutta2012Ex31}&Not observed&NA  & 22 & \texttt{DempeDutta2012Ex31}&0.525&2.92&3.26&137&115\\\rowcolor[gray]{.9}
23 & \texttt{DempeEtal2012}	&470&4.55e+09 & - & - & - & - & - & - & - \\\hline
24 & \texttt{DempeFranke2011Ex41}&492&1.33e+10 & 24 & \texttt{DempeFranke2011Ex41}&0.704&1.15&0.97&44.5&92\\\rowcolor[gray]{.9}
25 & \texttt{DempeFranke2011Ex42}&495&1.54e+10 & 25 & \texttt{DempeFranke2011Ex42}	&0.67&1.14&0.969&54.1&96\\\hline
26 & \texttt{DempeFranke2014Ex38}&495&1.54e+10 & 26 & \texttt{DempeFranke2014Ex38}&0.551&1.92&1.7&79.9&104\\\rowcolor[gray]{.9}
27 & \texttt{DempeLohse2011Ex31a}&502&2.17e+10 & 27 & \texttt{DempeLohse2011Ex31a}&0.525&3.45&5.08&6.02&51 
\\\hline
28 & \texttt{DempeLohse2011Ex31b}&510&3.2e+10 & - & - & - & - & - & - & - \\\rowcolor[gray]{.9}
29 & \texttt{DeSilva1978}&495&1.54e+10 & 29 & \texttt{DeSilva1978}&0.551&0.743&0.884&69&101\\\hline
30& \texttt{FalkLiu1995}&440&1.05e+09 & 30 & \texttt{FalkLiu1995}&$1.1*10^4$&36.3&30.3&$1.1*10^4$&205\\\rowcolor[gray]{.9}
31& \texttt{FloudasEtal2013}&505&2.51e+10 & 31 & \texttt{FloudasEtal2013}&0.525&1.7&1.65&183&121\\\hline
32 & \texttt{FloudasZlobec1998}&510&3.2e+10 & - & - & - & - & - & - & - \\\rowcolor[gray]{.9}
33 & \texttt{GumusFloudas2001Ex1}	&530&8.5e+10 & - & - & - & - & - & - & - \\\hline
34& \texttt{GumusFloudas2001Ex3}&510&3.2e+10 & - & - & - & - & - & - & - \\\rowcolor[gray]{.9}
35& \texttt{GumusFloudas2001Ex4}&502&2.17e+10 & - & - & - & - & - & - & - \\\hline
36& \texttt{GumusFloudas2001Ex5} &495&1.54e+10 & 36 & \texttt{GumusFloudas2001Ex5}&3.04&0.929&0.842&38.4&89\\\rowcolor[gray]{.9}
37& \texttt{HatzEtal2013}&Not observed&NA & 37 & \texttt{HatzEtal2013}&0.608&1.09&1.02&6.02&51  
\\\hline
38& \texttt{HendersonQuandt1958}&Not observed&NA & 38 & \texttt{HendersonQuandt1958}&$5.64*10^7$&2.34&2.33&$5.64*10^7$&380\\\rowcolor[gray]{.9}
39& \texttt{HenrionSurowiec2011}&Not observed&NA & 39 & \texttt{HenrionSurowiec2011}&0.551&1.12&1&6.02&51  
\\\hline
40& \texttt{IshizukaAiyoshi1992a} & 495&1.54e+10 & - & - & - & - & - & - & - \\\rowcolor[gray]{.9}
41 & \texttt{KleniatiAdjiman2014Ex3}&445&1.34e+09 & - & - & - & - & - & - & - \\\hline
42 & \texttt{KleniatiAdjiman2014Ex4}&485&9.46e+09 & 42 & \texttt{KleniatiAdjiman2014Ex4}	&0.739&0.617&0.531&42.4&91\\\rowcolor[gray]{.9}
43 & \texttt{LamparSagrat2017Ex23}	&Not observed&NA & 43 & \texttt{LamparSagrat2017Ex23}	&0.525&0.891&1.02&137&115\\\hline
44 & \texttt{LamparSagrat2017Ex31} &Not observed&NA & 44 & \texttt{LamparSagrat2017Ex31}&0.525&1.37&1.15&6.02&51  
 \\\rowcolor[gray]{.9}
45 & \texttt{LamparSagrat2017Ex32}	&Not observed&NA & 45 & \texttt{LamparSagrat2017Ex32}	&0.551&1&1.01&284&130\\\hline
46 & \texttt{LamparSagrat2017Ex33}	&495&1.54e+10 & 46 & \texttt{LamparSagrat2017Ex33}	&0.551&1.2&1&102&109\\\rowcolor[gray]{.9}
47 & \texttt{LamparSagrat2017Ex35}	&Not observed&NA & 47 & \texttt{LamparSagrat2017Ex35}	&0.525&1.65&1.53&298&131\\\hline
48 & \texttt{LucchettiEtal1987}&495&1.54e+10 & 48 & \texttt{LucchettiEtal1987}&0.525&1.03&1.43&6.02&51  
\\\rowcolor[gray]{.9}
49 & \texttt{LuDebSinha2016a}	&Not observed&NA & 49 & \texttt{LuDebSinha2016a}	&0.943&1.58&0.278&6.02& 51 
\\\hline
50 & \texttt{LuDebSinha2016b}	&Not observed&NA & 50 & \texttt{LuDebSinha2016b}	&0.67&0.263&0.255&6.02&51 
 \\\rowcolor[gray]{.9}
51 & \texttt{LuDebSinha2016c}&Not observed&NA & - & - & - & - & - & - & - \\\hline
52 & \texttt{LuDebSinha2016d}	 &890&3.61e+18 & - & - & - & - & - & - & - \\\rowcolor[gray]{.9}
53 & \texttt{LuDebSinha2016e}&900&5.88e+18 & - & - & - & - & - & - & - \\\hline
54 & \texttt{LuDebSinha2016f}&Not observed&NA & - & - & - & - & - & - & - \\\rowcolor[gray]{.9}
55 & \texttt{MacalHurter1997}&Not observed&NA & 55 & \texttt{MacalHurter1997}&0.551&18.1&18&6.02&51 
\\\hline
56 & \texttt{Mirrlees1999}	&Not observed&NA & 56 & \texttt{Mirrlees1999}	&0.579&0.118&0.118&6.02&51 
\\\rowcolor[gray]{.9}
57 & \texttt{MitsosBarton2006Ex38}	&398&1.36e+08 & 57 & \texttt{MitsosBarton2006Ex38}	&6.64&1.19e-5&9.98e-6&7.32&55\\\hline
58 & \texttt{MitsosBarton2006Ex39}	&400&1.5e+08 & - & - & - & - & - & - & - \\\rowcolor[gray]{.9}
59 & \texttt{MitsosBarton2006Ex310}&470&4.55e+09 & 59 & \texttt{MitsosBarton2006Ex310}	&56.8&0.165&0.141&56.8&97\\\hline
60 & \texttt{MitsosBarton2006Ex311}	&452&1.89e+09 & - & - & - & - & - & - & - \\\rowcolor[gray]{.9}
61 & \texttt{MitsosBarton2006Ex312}	&470&4.55e+09 & 61 & \texttt{MitsosBarton2006Ex312}	&0.99&0.108&0.125&6.02&51 
 \\\hline
62 & \texttt{MitsosBarton2006Ex313}	&505&2.51e+10 & 62 & \texttt{MitsosBarton2006Ex313}	&0.579&1.76&1.67&54.1&96\\\rowcolor[gray]{.9}
63 & \texttt{MitsosBarton2006Ex314}	&460&2.79e+09 & 63 & \texttt{MitsosBarton2006Ex314}	 &0.855&0.173&0.188&11.9&65\\\hline
64 & \texttt{MitsosBarton2006Ex315}	&470&4.55e+09 & 64 & \texttt{MitsosBarton2006Ex315}	&26&0.242&0.203&56.8&97\\\rowcolor[gray]{.9}
65 & \texttt{MitsosBarton2006Ex316}	&Not observed&NA & 65 & \texttt{MitsosBarton2006Ex316}	&1.33&5.49e-6&5.49e-6&6.02&51 
\\\hline
66 & \texttt{MitsosBarton2006Ex317}	&420&3.97e+08 & 66 & \texttt{MitsosBarton2006Ex317}	&3.7&0.00124&0.00107&6.02&51 
\\\rowcolor[gray]{.9}
67 & \texttt{MitsosBarton2006Ex318}	&Not observed&NA & 67 & \texttt{MitsosBarton2006Ex318}	&1.04&3.43e-9&3.43e-9&6.02&51 
\\\hline
68 & \texttt{MitsosBarton2006Ex319}	&398&1.36e+08 & - & - & - & - & - & - & - \\\rowcolor[gray]{.9}
69 & \texttt{MitsosBarton2006Ex320}	&485&9.46e+09 & 69 & \texttt{MitsosBarton2006Ex320}	&0.943&0.0434&0.0577&6.32&52\\\hline
70 & \texttt{MitsosBarton2006Ex321}&452&1.89e+09 & 70 & \texttt{MitsosBarton2006Ex321}	&1.09&0.0456&0.0381&10.3&62\\\rowcolor[gray]{.9}
71 & \texttt{MitsosBarton2006Ex322}	 &470&4.55e+09 & 71 & \texttt{MitsosBarton2006Ex322}	 &0.99&0.165&0.193&20.4&76\\\hline
72 & \texttt{MitsosBarton2006Ex323} &505&2.51e+10 & - & - & - & - & - & - & - \\\rowcolor[gray]{.9}
73 & \texttt{MitsosBarton2006Ex324}	&495&1.54e+10 & 73 & \texttt{MitsosBarton2006Ex324}	&0.855&1.17&1&6.02&51 
\\\hline
74 & \texttt{MitsosBarton2006Ex325}	 &505&2.51e+10 & - & - & - & - & - & - & - \\\rowcolor[gray]{.9}
75 & \texttt{MitsosBarton2006Ex326}	&505&2.51e+10 & - & - & - & - & - & - & - \\\hline
76 & \texttt{MitsosBarton2006Ex327}	&475&5.81e+09 & 76 & \texttt{MitsosBarton2006Ex327}	&1.26&0.104&0.105&7.32&55\\\rowcolor[gray]{.9}
77 & \texttt{MitsosBarton2006Ex328}	&510&3.2e+10 & - & - & - & - & - & - & - \\\hline
78 & \texttt{MorganPatrone2006a}&500&1.97e+10 & 78 & \texttt{MorganPatrone2006a}	&0.525&2.86&3.21&6.02&51 
 \\\rowcolor[gray]{.9}
79 & \texttt{MorganPatrone2006b}&505&2.51e+10 & 79 & \texttt{MorganPatrone2006b}&0.551&1.96&2.06&6.02&51 
 \\\hline
80 & \texttt{MorganPatrone2006c}	&470&4.55e+09 & 80 & \texttt{MorganPatrone2006c}&56.8&0.0894&0.206&56.8&97\\\rowcolor[gray]{.9}
81 & \texttt{MuuQuy2003Ex1}&Not observed&NA & - & - & - & - & - & - & - \\\hline
82 & \texttt{MuuQuy2003Ex2}	&Not observed&NA & - & - & - & - & - & - & - \\\rowcolor[gray]{.9}
83 & \texttt{NieEtal2017Ex34}&520&5.22e+10 & 83 & \texttt{NieEtal2017Ex34}&0.551&1.61&1.98&118&112\\\hline
84 & \texttt{NieEtal2017Ex52}	&Not observed&NA & - & - & - & - & - & - & - \\\rowcolor[gray]{.9}
85 & \texttt{NieEtal2017Ex54}&495&1.54e+10 & - & - & - & - & - & - & - \\\hline
86 & \texttt{NieEtal2017Ex57}	&850&5.13e+17 & - & - & - & - & - & - & - \\\rowcolor[gray]{.9}
87 & \texttt{NieEtal2017Ex58}	&780&1.69e+16 & - & - & - & - & - & - & - \\\hline
88 & \texttt{NieEtal2017Ex61}	&Not observed&NA & - & - & - & - & - & - & - \\\rowcolor[gray]{.9}
89 & \texttt{Outrata1990Ex1a}&505&2.51e+10 & 89 & \texttt{Outrata1990Ex1a}&0.608&1.39&1.31&192&122\\\hline
90 & \texttt{Outrata1990Ex1b}&510&3.2e+10 & 90 & \texttt{Outrata1990Ex1b}&0.551&2.42&2.3&223&125\\\rowcolor[gray]{.9}
91 & \texttt{Outrata1990Ex1c}&495&1.54e+10 & 91 & \texttt{Outrata1990Ex1c}&0.99&1.43&1.28&718&149\\\hline
92 & \texttt{Outrata1990Ex1d}	&495&1.54e+10 & - & - & - & - & - & - & - \\\rowcolor[gray]{.9}
93 & \texttt{Outrata1990Ex1e}&495&1.54e+10 & 93 & \texttt{Outrata1990Ex1e}	&0.855&1.53&1.29&873&153\\\hline
94 & \texttt{Outrata1990Ex2a}&520&5.22e+10 & 94 & \texttt{Outrata1990Ex2a}	&0.551&2.27&2.57&245&127\\\rowcolor[gray]{.9}
95 & \texttt{Outrata1990Ex2b}	&470&4.55e+09 & 95 & \texttt{Outrata1990Ex2b}&1.78&0.179&0.153&6.02&51 
 \\\hline
96 & \texttt{Outrata1990Ex2c}&510&3.2e+10 & 96 & \texttt{Outrata1990Ex2c}&0.551&2.09&2.64&72.5&102\\\rowcolor[gray]{.9}
97 & \texttt{Outrata1990Ex2d}&520&5.22e+10 & - & - & - & - & - & - & - \\\hline
98 & \texttt{Outrata1990Ex2e}	&470&4.55e+09 & 98 & \texttt{Outrata1990Ex2e}	&0.898&0.806&0.694&6.02&51 
\\\rowcolor[gray]{.9}
99 & \texttt{Outrata1993Ex31}&520&5.22e+10 & 99 & \texttt{Outrata1993Ex31}&0.551&2.1&2.27&144&116\\\hline
100 & \texttt{Outrata1993Ex32}	&680&1.28e+14 & - & - & - & - & - & - & - \\\rowcolor[gray]{.9}
101 & \texttt{Outrata1994Ex31}&910&9.58e+18 & - & - & - & - & - & - & - \\\hline
102 & \texttt{OutrataCervinka2009}&505&2.51e+10 & - & - & - & - & - & - & - \\\rowcolor[gray]{.9}
103 & \texttt{PaulaviciusEtal2017a}&400&1.5e+08 & 103 & \texttt{PaulaviciusEtal2017a}&107&1.58e-5&1.33e-5&107&110\\\hline
104 & \texttt{PaulaviciusEtal2017b}	&490&1.21e+10 & - & - & - & - & - & - & - \\\rowcolor[gray]{.9}
105 & \texttt{SahinCiric1998Ex2}&510&3.2e+10 & - & - & - & - & - & - & - \\\hline
106 & \texttt{ShimizuAiyoshi1981Ex1}&495&1.54e+10 & 106 & \texttt{ShimizuAiyoshi1981Ex1}&46.7&0.613&0.53&46.7&93\\\rowcolor[gray]{.9}
107 & \texttt{ShimizuAiyoshi1981Ex2}&520&5.22e+10 & - & - & - & - & - & - & - \\\hline
108 & \texttt{ShimizuEtal1997a}&910&9.58e+18 & 108 & \texttt{ShimizuEtal1997a}&0.638&4.2&3.53&212&124\\\rowcolor[gray]{.9}
109 & \texttt{ShimizuEtal1997b}&Not observed&NA & - & - & - & - & - & - & - \\\hline
110 & \texttt{SinhaMaloDeb2014TP3}&Not observed&NA & 110 & \texttt{SinhaMaloDeb2014TP3}&0.525&10.2&8.86&651&147\\\rowcolor[gray]{.9}
111 & \texttt{SinhaMaloDeb2014TP6}	&530&8.5e+10 & - & - & - & - & - & - & - \\\hline
112 & \texttt{SinhaMaloDeb2014TP7}&Not observed&NA & 112 & \texttt{SinhaMaloDeb2014TP7}&$5.75*10^{10}$&34.8&29&$5.75*10^{10}$&522\\\rowcolor[gray]{.9}
113 & \texttt{SinhaMaloDeb2014TP8}&520&5.22e+10 & - & - & - & - & - & - & - \\\hline
114 & \texttt{SinhaMaloDeb2014TP9}&505&2.51e+10 & - & - & - & - & - & - & - \\\rowcolor[gray]{.9}
115 & \texttt{SinhaMaloDeb2014TP10}	&480&7.41e+09 & - & - & - & - & - & - & - \\\hline
116 & \texttt{TuyEtal2007}	&495&1.54e+10 & 116 & \texttt{TuyEtal2007}&2.9&1.35&1.13&16.8&72\\\rowcolor[gray]{.9}
117 & \texttt{Vogel2002}&Not observed&NA & - & - & - & - & - & - & - \\\hline
118 & \texttt{WanWangLv2011}	&520&5.22e+10 & - & - & - & - & - & - & - \\\rowcolor[gray]{.9}
119 & \texttt{YeZhu2010Ex42}	&Not observed&NA & 119 & \texttt{YeZhu2010Ex42}	&0.814&0.73&0.656&6.02&51 
 \\\hline
120 & \texttt{YeZhu2010Ex43}	&Not observed&NA & 120 & \texttt{YeZhu2010Ex43}&4.49&0.482&0.402&6.02&51\\\rowcolor[gray]{.9}
121 & \texttt{Yezza1996Ex31}	&460&2.79e+09 & 121 & \texttt{Yezza1996Ex31}&223&0.0272&0.0565&223&125\\\hline
122 & \texttt{Yezza1996Ex41}&490&1.21e+10 & 122 & \texttt{Yezza1996Ex41}&0.608&0.609&0.638&46.7&93 \\\rowcolor[gray]{.9}
123 & \texttt{Zlobec2001a}&360&2.12e+07 & - & - & - & - & - & - & - \\\hline
124 & \texttt{Zlobec2001b}&495&1.54e+10 & - & - & - & - & - & - & -
\label{combinedtable}
\end{longtable}
\end{footnotesize}

We are now going to proceed with finding the threshold of $\lambda$ for which we start getting a solution and retain the value for larger values of $\lambda$. We are going to proceed with the technique of finding small threshold $\bar \lambda$ and large threshold $\lambda^*$, which was briefly discussed earlier.
To find the threshold  we first extract the value of the final $Error^*$ for each example from \cite{LMnumer20}.
To find $\bar \lambda$ we run the algorithm with $\lambda$ being defined as $\lambda:=0.5\times1.05^k$, with the new stopping criteria: 
$$\textbf{Stop if } \phantom{-} Error\leq1.1 Error^*.$$
Of course, we also relax all of the aforementioned stopping criteria, as $Error$ is the main measure here and we know that desired value of $Error$ exists.  We then define $\bar \lambda := 0.5 \times 1.05^{\bar k}$, where $\bar k$ is the number of iterations after stopping whenever we get $Error<1.1Error^*$.
For most of the cases we detect $\bar \lambda$ early due to a good solution after the first few iterations in the same manner as shown in Figure  \ref{twofigsexample2} (c).
Since for many examples we satisfy condition  $Error\leq1.1 Error^*$ for early iterations (before the inflection point is achieved),  we further introduce $\lambda^*$, \emph{the large threshold} $\lambda$. The value of $\lambda^*$ will be used to represent the value of the penalty parameter at the inflection point, where solution starts to be recovered for large values of $\lambda$. This will be obtained in the same way as $\bar \lambda$ with the only difference that we additionally impose the condition to stop after at least 50 iterations. To obtain $\lambda^*$ we run the algorithm with $\lambda:=0.5\times1.05^k$ and the following stopping criteria:
$$\textbf{Stop if } \phantom{-} Error\leq1.1 Error^*  \phantom{-} \&   \phantom{-}iter>50.$$
Then the large threshold is defined as $\lambda^* := 0.5\times 1.05^{k^*}$, where $k^*$ is the number of iterations after stopping whenever we get $Error<1.1Error^*$ and $k>50$.
This way $\bar \lambda$ represents the first (smallest) value of $\lambda$ for which good solution was obtained, while $\lambda^*$ represent the actual threshold after which solution is retained for large values of $\lambda$.
The demonstration of stopping at the inflection point for large threshold $\lambda^*$ was shown in Figure \ref{twofigsexample2} (b) and for small threshold $\bar \lambda$ in Figure  \ref{twofigsexample2} (d).

It makes sense to test only the examples where algorithm performed well and produced a good solution. For the rest of the examples, the value of $Error^*$ would not make sense as the measure to stop, and we do not obtain good solutions for these examples by the algorithm anyway.
From the optimistic perspective we could consider recovered solutions to be the solutions for which the optimal value of upper-level objective was recovered with some prescribed tolerance.
Taking the tolerance of $20\%$, the total amount of recovered solutions by the method with varying $\lambda$ is 72/117  ($61.54\%$) for the cases where solution was reported in BOLIB \cite{bolib17}. This result will be shown in more details in Section \ref{UpperlevSec}.
Let us look at the thresholds $\bar \lambda$ and $\lambda^*$ for these examples in the figure below, where the value of $\lambda$ is shown on $y-axis$ and the example numbers on the $x-axis$.
\begin{figure}[H]
\centering
\includegraphics[scale=0.5]{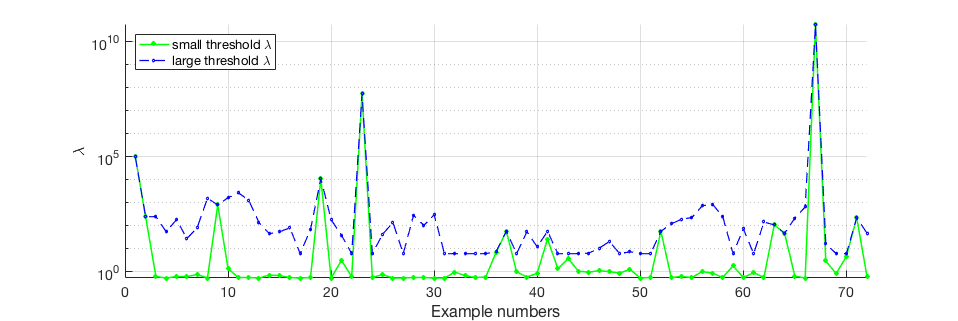}
    \caption{Small threshold $\bar \lambda$ and large threshold $\lambda^*$ for examples with good solutions}%
    \label{thresholdslamfig}%
\end{figure}
For $68/72$ problems, we observe that large threshold of the penalty parameter has the value $\lambda^* \leq 10^4$, which shows that we usually can obtain a good solution before 205 iterations. This further justifies stopping algorithm after 200 iterations (see next section) if there is no significant step to step improvement.
As for the main outcome of Figure \ref{thresholdslamfig}, we observe that small threshold is smaller than large threshold ($\bar \lambda < \lambda^*$) for $59/72$ ($82 \%$) problems.
This clearly shows that for majority of the problems, for which we recover solution with $\lambda:=0.5\times1.05^k$, we obtain a good solution for small $\lambda$ as well as for large $\lambda$. This demonstrates that small $\lambda$ could in principle be good for the method.
For the rest $18\%$ of the problems we have $\bar \lambda = \lambda^*$, meaning that good solution was not obtained for $\lambda<6$ for these examples.
 This also means that we typically obtain good solution for small values of $\lambda$ and for large values of $\lambda$, but not for the medium values ($\bar \lambda<\lambda<\lambda^*$).

As for the general observations of Figure \ref{thresholdslamfig}, for $42/72$ ($58.33\%$) examples we observed that the large threshold $\lambda^*$ is located somewhere in between $90-176$ iterations with $40<\lambda^*<2680$.
For $7/72$ problems threshold is in the range $6.02<\lambda^* \leq40$, and for only $4/72$ problems $\lambda^*>1.1\times10^4$.
Once again, this justifies that typically $\lambda$ does not need to be large.
It also suggests the optimal values of $\lambda$ for the tested examples, at least for our solution method.
We observe that for $19/72$ problems we get $\lambda^*=51$, which should be treated carefully as this could mean that the inflection point could possibly be achieved before $50$ for these examples. Nevertheless, $\lambda^*=6.02$ is still a good value of penalty parameter for these examples as solutions are retained for $\lambda>\lambda^*$.

As going to be observed in Section \ref{SecNumExper} we could actually argue that smaller values of $\lambda$ work better for our method not only for varying $\lambda$ but also for fixed $\lambda$. Together with the fact that we often have the behaviour as demonstrated in Figure \ref{twofigsexample2} (c), it follows that small $\lambda$ could be more attractive for the method we implement. We even get better values of $Error$ and better solutions for small values of $\lambda$ for some examples. Hence we draw the conclusion that small values of $\lambda$ can generate good solutions.
Since it is typical to use large values of $\lambda$ for other penalization methods  (e.g. in  \cite{bert82,burke91,pillo89}), it is interesting what could be the reasons that small $\lambda$ worked better for our case.
This could be due to the specific nature of the method, or due to the fact that we do not do full penalization in the usual sense. Other reason could come from the structure of the problems in the test set. The exact reason of why such behaviour was observed remains an open question. What is important is that this could possibly be the case that small values of $\lambda$ would be good for some other penalty methods and optimization problems of different nature.
This result contradicts typical choice of large penalty parameter for general penalization methods for optimization problems.
 As the conclusion for our framework, we can claim that for our method $\lambda$ needs not to be large.

\subsection{Partially calm examples}
\begin{figure}[H]
   \hspace{-0.44cm}
    \subfloat[Example 8]{{\includegraphics[width=6cm]{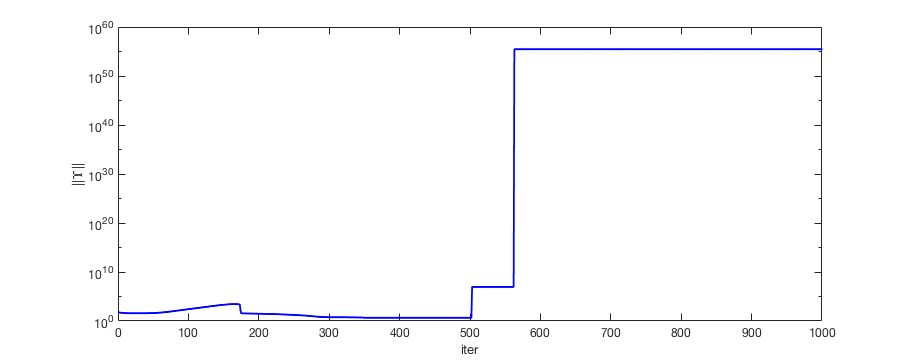} }}%
   \subfloat[Example 123]{{\includegraphics[width=6.4cm]{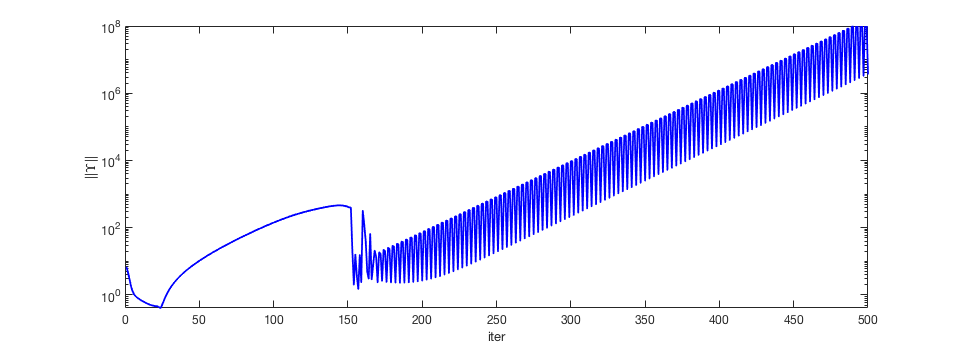} }}%
   \subfloat[Example 45]{{\includegraphics[width=5.2cm]{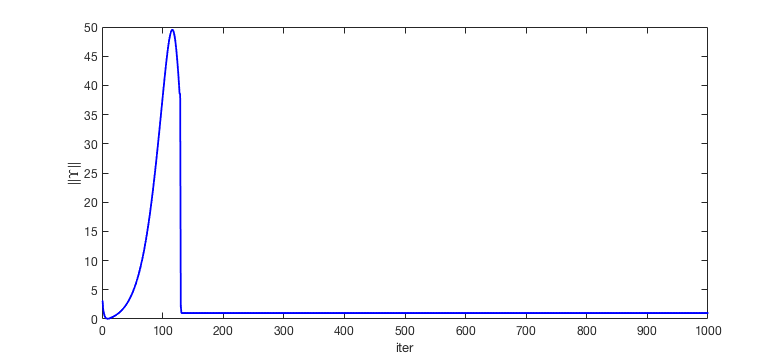} }}%
    \caption{(a) and (b) are obtained for 1000 iterations, while (c) is based on 500 iterations.}%
    \label{fig:awesomeImage3}%
\end{figure}

Intuitively, one would think that for partially calm examples, Algorithm \ref{algorithm LMLs} would behave well, in the sense that varying $\lambda$ increasingly would lead to a good convergence behaviour. To show that it is not necessarily the case, we start by considering the following result identifying a class of bilevel program of the form \eqref{initialbilev} that is automatically partially calm.
\begin{thm}[\cite{MMZ20}] \label{CalmExamples} Consider a bilevel program \eqref{initialbilev}, where $G$ is independent of $y$ and the lower-level optimal solution map is defined as follows, with $c\in \mathbb{R}^m$, $d\in \mathbb{R}^p$, $B\in \mathbb{R}^{p\times m}$, and $A :\mathbb{R}^n \rightarrow \mathbb{R}^p$:
\begin{equation}\label{C2}
S(x):= \arg \min_y \left\{c^T y |\; A(x)+B y\leq d\right\}.
\end{equation}
In this case, problem \eqref{initialbilev} is partially calm at any of its local optimal solutions.
\end{thm}
Examples 8, 40, 43, 45, 46, 188, and 123 in the BOLIB library (see Table \ref{combinedtable}) are of the form described in this result.
The expectation is that these examples will follow the pattern of retaining solution after some threshold, that is for  $\lambda>\lambda^*$, as they fit the theoretical structure behind the penalty approach as described in Theorem \ref{thresholdlem}.
Note that all of these examples follow the pattern shown in Figure \ref{alpha2}(a).
However, if we relax the stopping criteria used to mitigate the effects of ill-conditioning, as discussed in the previous two subsections,  varying $\lambda$  for 1000 iterations for these seven partially calm examples leads to the 3 typical scenarios demonstrated  in Figure \ref{fig:awesomeImage3}. 

In the first case of Figure \ref{fig:awesomeImage3}, we can clearly see the algorithm is performing well, retaining the solution for the number of iteration, but then blows up at one point (after 500 iterations) and never goes back to reasonable solution values. Examples 40 and 46 also follow this pattern.
Example 123 (second picture in Figure \ref{fig:awesomeImage3}) shows a slightly different picture, where the zig-zagging pattern is observed.
Algorithm \ref{algorithm LMLs} blows up at some point and starts zig-zagging away from the solution after obtaining it for a smaller value of $\lambda$. Zig-zagging is very common issue in penalty methods and often caused by ill-conditioning \cite{nocedal99}. Note that Example 118 exhibits a similar behaviour. This is somewhat similar to scenario 1. However, we put this separately as zig-zagging issue is often referred to as the danger that could be caused by ill-conditioning of a penalty function.
The last picture of Figure \ref{fig:awesomeImage3} shows a case where Algorithm \ref{algorithm LMLs} runs very well without any ill-behaviour observed for all the 1000 iterations. It could be possible that the algorithm could blow up after more iterations if we keep increasing $\lambda$.
It could also be possible that ill-conditioning  does not occur for this example at all, as the Hessian of $\Upsilon^\lambda$ \eqref{bilevncp21} is not affected by $\lambda$. Out of the seven BOLIB problems considered here, only examples 43 and 45 follows this pattern.

\section[Performance comparison under fixed and varying penalty parameter]{Performance comparison for the Levenberg-Marquardt method under fixed and varying partial exact penalty parameter}\label{SecNumExper}
Following the discussion from the previous section, two approaches for selecting the penalty parameter $\lambda$ are tested and compared on the Levenberg-Marquardt method.
Recall that for the varying $\lambda$ case, we define the penalty parameter sequence as $\lambda_k := 0.5\times1.05^k$, where $k$ is the iteration number.
When fixed values of penalty parameter are considered, ten different values of $\lambda$ are used for the experiments; i.e., $\lambda \in \{10^6, 10^5, ..., 10^{-3}\}$. 
For the fixed values of $\lambda$ one could choose best $\lambda$ for each example to see if at least one of the selected values worked well to recover the solution.
As in the previous section, the examples used for the experiments are from the Bilvel Optimization LIBrary of Test Problems (BOLIB) \cite{bolib17}, which contains 124 nonlinear examples. The experiments are run in MATLAB, version R2016b, on MACI64. Here, we present a summary of the results obtained; more details for each example are reported in the supplementary material \cite{LMnumer20}.  It is important to mention that algorithm always converges, unlike its Gauss-Newton counterpart studied in \cite{FTZ20}, where the corresponding algorithm diverges for some values of $\lambda$.

\subsection{Practical implementation details}
For Step 0 of Algorithm \ref{algorithm LMLs} we set the tolerance to $\epsilon:=10^{-5}$ and the maximum number of iterations to be $K:=1000$. We also choose $\alpha_0 := \norm{\Upsilon^\lambda(z^0)}$ (if $\lambda$ is varying, we set $\lambda:=\lambda_0$ here), $\gamma_0:= 1$, $\rho=0.5$, and  $\sigma=10^{-2}$. The selection of $\sigma$ is based on the overall performance of the algorithm while the other parameters are standard in the literature.
${}$\\

\noindent \textbf{\emph{Starting point.}} The experiments have shown that the algorithm performs much better if the starting point $(x^0,y^0)$ is feasible. As a default setup, we start with $x^0=1_n$ and $y^0=1_m$. If the default starting point does not satisfy at least one constraint and algorithm diverges, we choose a feasible starting point; see \cite{LMnumer20} for such choices. To be more precise, if for some $i$, $G_i(x^0,y^0)>0$ or for some $j$ we have $g_j(x^0,y^0)>0$ and the algorithm does not converge to a reasonable solution, we generate a starting point such that  $G_i(x^0,y^0)=0$ or $g_j (x^0,y^0) =0$. Subsequently, the Lagrange multipliers are initialised at $u^0 = \max\,\left\{0.01,\;-g(x^0,y^0)\right\}$,  $v^0=\max\,\left\{0.01,\;-G(x,y)\right\}$, and $u^0 = w^0$.
${}$\\

\noindent \textbf{\emph{Smoothing parameter $\mu$.}}
The smoothing process, i.e., the use of the parameter $\mu$, is only applied when necessary (Step 2 and Step 3), where the derivative evaluation for $\Upsilon^\lambda_\mu$ is required. We tried both the case where $\mu$ is fixed to be a small constant for all iterations, and the situation where the smoothing parameter is sequence $\mu_k \downarrow 0$. 
Setting the decreasing sequence to $\mu_k :=0.001/(1.5^k)$, and testing its behaviour and comparing it with a fixed small value ($\mu:=10^{-11}$), in the context of Algorithm \ref{algorithm LMLs}, we observed that both options lead to almost the same results. 
Hence, we stick to what is theoretically more suitable, that is the smoothing decreasing sequence ; i.e., $\mu_k:=0.001/(1.5^k)$.
${}$\\

\noindent \textbf{\emph{Descent direction check and update.}} For the sake of efficiency, we calculate the direction $d^k$ by solving \eqref{SLMdirnum} with the Gaussian elimination. Considering the line search in Step 3, if we have $\norm{\Upsilon^{\lambda}(z^k+\gamma_k d^k)}^2 < \norm{\Upsilon^{\lambda}(z^k)}^2 + \sigma \gamma_k \nabla \Upsilon^\lambda(z^k)^T \Upsilon^\lambda(z^k) d^k$, we redefine $\gamma_k = \gamma_k / 2$ and check again.
Recall that the Levenberg-Marquardt direction can be interpreted as a combination of Gauss-Newton and steepest descent directions.
In fact, if $\alpha_k = 0$ this direction is a Gauss-Newton direction one and when as $\alpha_k \to \infty$ the direction $d^k$ from  \eqref{SLMdirnum} tends to a steepest descent direction.
Hence, if the Levenberg-Marquardt direction is not a descending at some iteration, we give more weight to the steepest descent direction.
Hence, when $\norm{\Upsilon^\lambda(z^k)}>\norm{\Upsilon^\lambda(z^{k-1})}$
setting $\alpha_{k+1} := 10000 \norm{\Upsilon^\lambda(z^k)}$ has led to an oeverall good performance of Algorithm \ref{algorithm LMLs} for test set used in this paper.
${}$\\

\noindent \textbf{\emph{Stopping Criteria.}}  The primary stopping criterion for Algorithm \ref{algorithm LMLs} is $\norm{\Upsilon^{\lambda}_\mu (z^k)}<\epsilon$, as requested in Step 1. However, robust safeguards are needed to deal with ill-behaviours typically due to the size of the penalty parameter $\lambda$. Hence, for the practical implementation of the method, we  set $\epsilon=10^{-5}$ and \textbf{stop} if one of the following six conditions is satisfied:
\begin{enumerate}
\item $\norm{\Upsilon^\lambda(z^k)}<\epsilon$,
\item $\big| \norm{\Upsilon^\lambda(z^{k-1})} - \norm{\Upsilon^\lambda(z^{k})} \big|<10^{-9}$,
\item $\big| \norm{\Upsilon^\lambda(z^{k-1})} - \norm{\Upsilon^\lambda(z^{k})} \big|<10^{-4}$ and  $iter>200$,
\item $\norm{\Upsilon^\lambda(z^{k-1})} - \norm{\Upsilon^\lambda(z^{k})} <0$ and $\norm{\Upsilon^\lambda(z^k)}< 10$ and $iter>175$,
 \item $\norm{\Upsilon^\lambda(z^k)}<10^{-2}$ and $iter>500$,
\item  $\norm{\Upsilon^\lambda(z^k)}>10^2$ and $iter>200$.
\end{enumerate}
The additional stopping criteria is important to ensure that algorithm is not running for too long. The danger of running algorithm for too long is that ill-conditioning could occur. Further, we typically observe the pattern that we recover solution earlier than algorithm stops. This appears due to the nature of the overdetermined system. We do not know beforehand the tolerance with which we can solve for $\norm{\Upsilon^\lambda(z)}$ as $\Upsilon^\lambda$ is overdetermined system. Hence, it is hard to select $\epsilon$ that would fit all examples and allow to solve examples with efficient tolerance. With the stopping criteria defined above we avoid running unnecessary iterations, retaining the obtained solution.
To avoid algorithm running for too long and to prevent $\lambda$ to become too large, we impose additional stopping criterion 3, 5 or 6 above.
These criteria are motivated by the behaviour of the algorithm.

For almost all of the examples we observe that after 100-150 iterations we obtain the value reasonably close to the solution but we cannot know beforehand what would be the tolerance of $\norm{\Upsilon^\lambda(z^k)}$ to stop. Choosing small $\epsilon$ would not always work due to the overdetermined nature of the system being solved. Choosing $\epsilon$ too big would lead to worse solutions and possibly not recover some of the solutions. Further, a quick check has shown that ill-conditioning issue typically takes place after 500 iterations for majority of the problems.
 For these reasons we stop if the improve of the $Error$ value from step to step becomes too small, $\big| \norm{\Upsilon^\lambda(z^{k-1})} - \norm{\Upsilon^\lambda(z^{k})} \big|<10^{-4}$, after the algorithm has performed $200$ iterations. Since ill-conditioning is likely to happen after $500$ iterations we stop if by that time we obtain a reasonably small $Error$, $\norm{\Upsilon^\lambda(z^k)}<10^{-2} $. Finally, if it turns out that system cannot be solved with a good tolerance, such that we would obtain a reasonably small value of the $Error$, we stop if the Error after $200$ iterations is big, $\norm{\Upsilon^\lambda(z^k)}>10^2$. This way additional stopping criteria plays the role of safeguard to prevent ill-conditioning and also does not allow the algorithm to keep running for too long once a good solution is obtained.

\subsection{Accuracy of the upper-level objective function}\label{UpperlevSec}
Here, we compare the values of the upper-level objective functions at points computed by the Levenberg-Marquardt algorithm with fixed $\lambda$ and varying $\lambda$.
 For the comparison purpose, we focus our attention only on  117 BOLIB examples \cite{bolib17}, as solutions are not known for the other seven problems.
 To proceed, let $\bar{F}_A$ be the value of upper-level objective function at the point $(\bar{x},\bar{y})$ obtained by the algorithm and $\bar{F}_K$  the value of this function at the known best solution point reported in the literature (see corresponding references in \cite{bolib17}).
We consider all fixed $\lambda\in\{10^6,10^5,...,10^{-3}\}$ and  varying $\lambda$ in one graph and present the results in Figure \ref{upperlevvarandfixedlam} below, where we have the relative error $(\bar{F}_A - \bar{F}_K) / (1+|\bar{F}_K|)$ on the $y-axis$ and number of examples on the x-axis, starting from $30th$ example. We further plot the results for the best fixed value of $\lambda$. The graph is plotted in the order of increasing error.
From the Figure \ref{upperlevvarandfixedlam} above we can clearly see that much more solutions were recovered for the small values of fixed $\lambda$ than for large values. For instance with the allowable accuracy error of $\leq20\%$ we recover solutions for $78.63\%$ for fixed $\lambda\in\{10^{-2},10^{-3}\}$, while for $\lambda\in\{10^6,10^{5},10^4,10^3,10^2\}$ we recover at most $40.17\%$ solutions. Interestingly, the worst performance is observed for fixed $\lambda=100$.
With the varying $\lambda$ we observe that algorithm performed averagely in comparison between large and small fixed values of $\lambda$, recovering $59.83\%$ of the solutions with the accuracy error of $\leq20\%$.
It is worth saying that implementing Algorithm \ref{algorithm LMLs} with $\lambda:=0.5\times 1.05^k$ still recovers over half of the solutions, which is not too bad. However, fixing $\lambda$ to be small recovers way more solutions, which shows that varying $\lambda$ is not the most efficient option for our case.

It was further observed that for some examples only $\lambda \geq 10^3$ performed well, while for others small values ($\lambda<1$) showed good performance.
If we were able to pick best fixed $\lambda$ for each example, we would obtain negligible (less than $10 \%$) error for upper-level objective function for $85.47\%$ of the tested problems. With the accuracy error of $\leq25\%$ our algorithm recovered solutions for $88.9\%$ of the problems for the best fixed $\lambda$ and for $61.54\%$ with varying $\lambda$. This means that if one can choose the best fixed $\lambda$ from the set of different values, fixing $\lambda$ is much more attractive for the algorithm. It was further observed that for some examples only $\lambda \geq 10^3$ performed well, while for others small values of $\lambda$ ($\lambda<1$) showed good performance. However, if one does not have a way to choose the best value or a set of potential values cannot be constructed efficiently for certain problems, varying $\lambda$ could be a better option to choose. Nevertheless, for the test set of small problems from BOLIB \cite{bolib17}, fixing $\lambda$ to be small performed much better than varying $\lambda$ as increasing sequence. Further, if one could run the algorithm for all fixed $\lambda$ and was able to choose the best one, the algorithm with fixed $\lambda$ performs extremely well compared to varying $\lambda$. In other words, algorithm almost always finds a good solution for at least one value of fixed $\lambda \in \{10^6, 10^5, ..., 10^{-3}\}$.

\begin{figure}[H]
\hspace{+0.1cm}
\includegraphics[width=17cm]{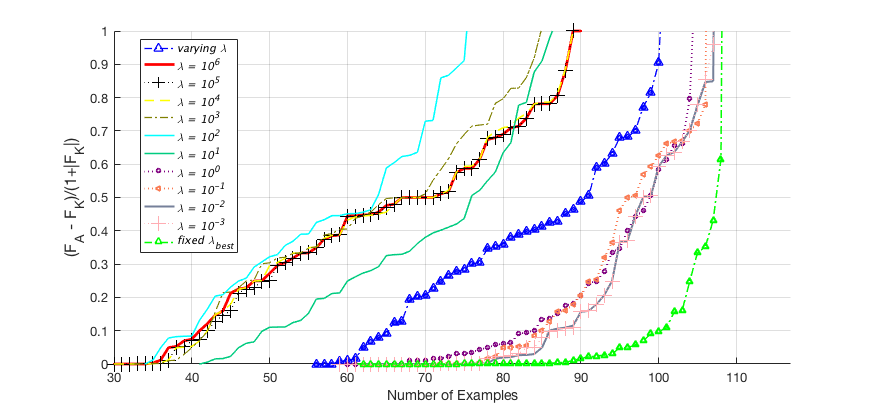} %
\caption{Error of the upper-level objective value for examples with known solutions}%
    \label{upperlevvarandfixedlam}%
\end{figure}

\subsection{Feasibility check}

\begin{figure}[H]
\hspace{-0.1cm}
\includegraphics[width=17cm]{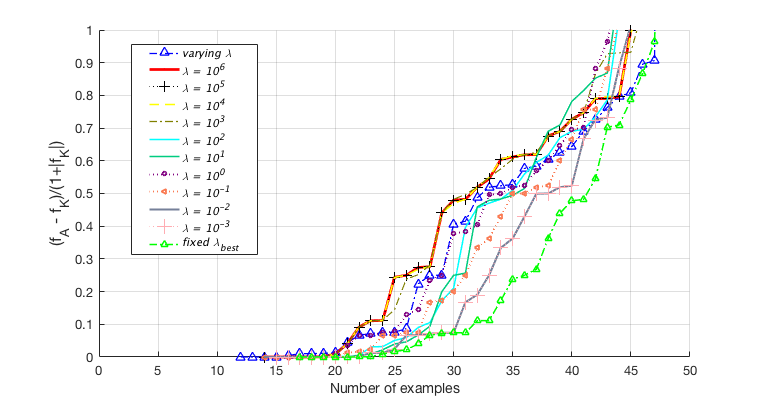} %
 \caption{Feasibility error for the lower-level problem in increasing order}%
    \label{lowerlevvarandfixedlam}%
\end{figure}

Considering the structure of the feasible set of problem \eqref{initialvalfuncform0}, it is critical to check whether the points computed by our algorithms satisfy the value function constraint $f(x,y)\leq \varphi(x)$, as it is not explicitly included in the expression of $\Upsilon^\lambda$ \eqref{bilevncp21}.
If the lower-level problem is convex in $y$ and a solution generated by our algorithms satisfies \eqref{kktbilev32} and \eqref{kktbilev62}, then it will verify the value function constraint. Conversely, to guaranty that a point $(x, y)$ such that $y\in S(x)$ satisfies \eqref{kktbilev32} and \eqref{kktbilev62}, a constraint qualification (CQ) is necessary. Note that conditions \eqref{kktbilev32} and \eqref{kktbilev62} are incorporated in the stopping criterion of Algorithm \ref{algorithm LMLs}. To check whether the points obtained are feasible, we first identify the BOLIB examples, where the lower-level problem is convex w.r.t. $y$.
As shown in \cite{FTZ20} it turns out that a significant number of test examples have linear lower-level constraints. For these examples, the lower-level convexity is automatically satisfied.
We detect 49 examples for which some of these assumptions are not satisfied, that is the problems having non-convex lower-level objective or some of the lower-level constraints being nonconvex. For these examples, we compare the obtained solutions  with the known ones from the literature.
Let $f_A$ stand for $f(\bar{x},\bar{y})$ obtained by one of the tested algorithms and $f_K$ to be the known optimal value of lower-level objective function.
In the graph below we have the lower-level relative error, $(f_A - f_K) / (1+|f_K|)$, on the y-axis, where the error is plotted in increasing order.
In Figure \ref{lowerlevvarandfixedlam} below we present results for all fixed $\lambda\in\{10^6,10^5,...,10^{-3}\}$ as well as varying $\lambda$ defined as $\lambda:=0.5\times1.05^k$.

From the Figure \ref{lowerlevvarandfixedlam} above we can see that for 20 problems the relative error of lower-level objective is negligible ($<5\%)$ for all values of fixed $\lambda$ and varying $\lambda$.
We have seen that convexity and a CQ hold for the lower-level hold for 74 test examples.
We consider solutions for these problems to be feasible for the lower-level problem.
Taking satisfying feasibility error to be $<20 \%$  and using information from the graph above, we claim that feasibility is satisfied for at most $100$ ($80.65\%$) problems for fixed $\lambda\in \{10^6,10^5,10^4,10^3\}$, for $101-104$ ($81.45-83.87\%$) problems for $\lambda\in \{10^3,10^2,10^1,10^0,10^{-1}\}$ and for $106$ ($85.48\%$) problems for $\lambda\in \{10^{-2},10^{-3}\}$. We further observe that feasibility is satisfied for 101 (81.4\%) problems for varying $\lambda$.
Considering we could choose best fixed $\lambda$ for each of the examples, we could also claim that feasibility is satisfied for 108 (87.1\%) problems for best fixed $\lambda$.
From Figure \ref{lowerlevvarandfixedlam} we note that slightly better feasibility was observed for smaller values of fixed $\lambda$ than for the big ones and that varying $\lambda$ has shown average performance between these magnitudes in terms of the feasibility.

\subsection{Experimental order of convergence}
Recall that the experimental order of convergence (EOC) is defined by
$$
\mbox{EOC}:=\max\left\{\frac{\log\|\Upsilon^\lambda(z^{K-1})\|}{\log\|\Upsilon^\lambda(z^{K-2})\|},\, \frac{\log\|\Upsilon^\lambda(z^{K})\|}{\log\|\Upsilon^\lambda(z^{K-1})\|} \right\},
$$
where $K$ is the number of the last iteration \cite{newtonbilevel18}.
If $K=1$, no EOC will be calculated (EOC$=\infty$).
EOC is important to estimate the local behaviour of the algorithm and to show whether this practical convergence reflects the theoretical convergence result stated earlier.
Let us  consider EOC for fixed $\lambda\in \{10^6,...,10^{-3}\}$ and for varying $\lambda$ ($\lambda=0.5\times1.05^k$) in Figure  \ref{EOCvarandfixedlam} below.




It is clear from this picture that for most of the examples our method has shown linear experimental convergence.
This is slightly below the quadratic convergence established by Theorem \ref{LMconvtheoremEBs1}. It is however important to note that the method always converges to a number, although sometimes the output might not be the optimal point for the problem.
There are a few examples that shown better convergence for each value of $\lambda$, with the best ones being $\lambda\in\{10^{-3},10^{-2},10^{-1},10^6\}$ as seen in the figure above.
These fixed values have shown slightly better EOC performance than varying $\lambda$. Varying $\lambda$ showed slightly better convergence than fixed $\lambda\in\{10^0,10^1,10^2,10^3,10^4,10^5\}$.
EOC bigger than $1.2$ has been obtained for less than 5 (4.03 \%) examples for fixed $\lambda\in\{10^0,10^1,10^2,10^3,10^4,10^5\}$, while varying $\lambda$ showed such EOC for 11 (8.87\%) examples. Fixed $\lambda=10^6$ has shown almost the same result as varying $\lambda$ with rate of convergence greater than $1.2$ for 12  (9.67\%) examples, while $\lambda=10^{-1}$ has demonstrated such EOC for 14 (11.29\%) examples and $\lambda\in\{10^{-2},10^{-3}\}$ for 17 (13.71 \%) examples.
Finally, in the graph above, we can see that for all values of $\lambda$ only a few ($\leq4/124$) examples have worse than linear convergence.

\begin{figure}[H]
\hspace{-0.4cm}
\includegraphics[width=18cm]{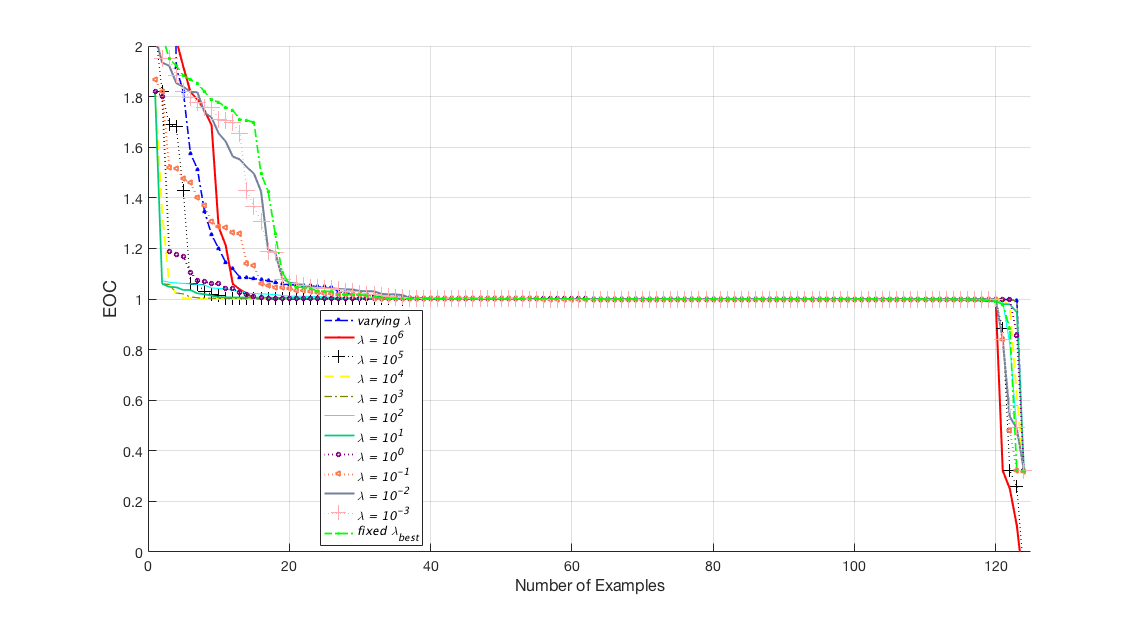} %
 \caption{Observed EOC at the last iterations for all examples (in decreasing order)}%
    \label{EOCvarandfixedlam}%
\end{figure}

\subsection{Line search stepsize}
Let us now look at the line search stepsize, $\gamma_k$, at the last step of the algorithm for each example.
Consider all fixed $\lambda$ and varying $\lambda$  in Figure \ref{Stepvarandfixedlam} below.
This is quite important to know two things.
Firstly, how often line search was used at the last iteration, that is how often implementation of line search was clearly important.
Secondly, as main convergence results are for the pure method this would be demonstrative to note how often the pure (full) step was made at the last iteration. This can then be compared with the experimental convergence results in the previous subsection, namely with Figure \ref{EOCvarandfixedlam}.
\begin{figure}[H]
\hspace{-0.2cm}
\includegraphics[width=17.5cm]{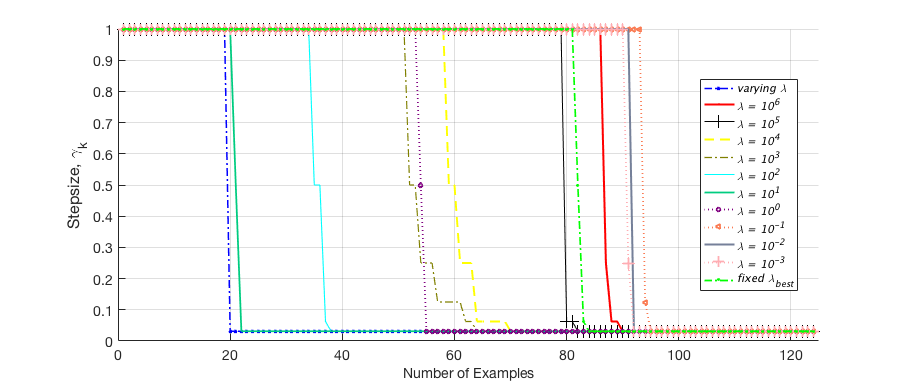} %
    \caption{Stepsize made at the last iteration for all examples (in decreasing order)}%
    \label{Stepvarandfixedlam}%
\end{figure}
In the figure above whenever stepsize on the y-axis is equal to $1$ it means the full step was made at the last iteration. For these cases the convergence results shown in Theorem \ref{LMconvtheoremEBs1} could be considered valid.
From the graphs above we observe that stepsize at the last iteration was rather $\gamma_k=1$ or $\gamma_k<0.05$.
We observe that for varying $\lambda$ algorithm would typically do a small step at the last iteration. It seems that algorithm with varying $\lambda$ benefits more from line search technique than the algorithm with fixed $\lambda$. Possibly, pure Levenberg-Marquardt method with varying $\lambda$ would not converge for most of the problems.
 Interestingly, for fixed values of $\lambda$ stepsize was $\gamma_k<0.05$ at the last iteration much more often for the values of $\lambda$ that showed worse performance in terms of recovering solutions (i.e. $\lambda\in \{10^1, 10^2\}$).
We also observe that for medium values of $\lambda$ ($\lambda\in\{10^4,10^3,10^2,10^1,10^0\}$) full stepsize was made for less than half of the examples. For large values $\lambda\in\{10^5,10^5\}$  full step was made for $63.71\%$ and  $70.16\%$ of the problems respectively.
Further on, small values of $\lambda$ for which more solutions were recovered would do the full step at the last iteration for most of the examples. For instance, with $\lambda=10^{-3}$ and $\lambda=10^{-2}$ full step was made at the last iteration for $73.39\%$ of the problems, while for $\lambda=10^{-1}$ full step was made for $75.81\%$ of the problems.
In terms of the fixed $\lambda_{best}$ it is interesting to observe that full step was used only for 82/124 ($66.13 \%$) of the problems, meaning that for a third of the problems linesearch was implemented in the last step for the best tested value of $\lambda$.
This also coincides with the results of Figure  \ref{EOCvarandfixedlam} where with smaller values of $\lambda$ algorithm has shown faster than linear convergence for more examples than for big values of $\lambda$.
This is likely to be the case that small steps were made in the other instances due to the non-efficient direction of the method at the last iteration.


\section{Final comments}

We introduce a smoothing Levenberg-Marquardt method to solve LLVF-based optimality conditions for optimistic bilevel programs. Since these optimality conditions are parameterized by a certain $\lambda$, its selection is carefully studied in light of the BOLIB library test set.
 Surprisingly, small values of this partial exact penalty parameter showed a good behaviour for our method, as they generally performed very well, although classical literature on exact penalization usually suggest  large values.
However, relatively medium values of $\lambda$ did not perform as well.
Furthermore, as both the varying and fixed scenarios were considered for $\lambda$, it was observed that both approaches showed a linear experimental order of convergence for most of the examples.
 Average CPU time for all fixed $\lambda$ is $0.243$ seconds, while it is $0.193$ seconds for $\lambda_{best}$ and $0.525$ seconds for varying $\lambda$.
The algorithm with varying $\lambda$ turns out to be more than twice slower than for fixed $\lambda$. However, running it for all fixed values of $\lambda$ can take way more time than the case where the value of $\lambda$ is varying.

\appendix


\begin{thebibliography}{99}
\bibitem{Allendestill12}
G.B. Allende  and G. Still,
\newblock {\it Solving bi-level programs with the KKT-approach},
\newblock Mathematical Programming 131:37-48 (2012)



\bibitem{behling16}
R. Behling, A. Fischer, G. Haeser, A.Ramos and K. Sch\"{o}nefeld,
\newblock {\it On the constrained error bound condition and the projected Levenberg-Marquardt method},
\newblock Optimization 66(8): 1397-1411
(2016)

\bibitem{bert82}
D.P. Bertsekas,
\newblock Constrained optimization and Lagrange multiplier methods,
\newblock  Academic Press (1982)

\bibitem{burke91}
J.V. Burke,
\newblock {\it An exact penalization viewpoint of constrained optimization},
\newblock SIAM journal on control and optimization 29(4): 968-998
(1991)

\bibitem{bilevelmpec10}
S. Dempe  and J. Dutta,
\newblock {\it Is bilevel programming a special case of mathematical programming with equilibrium constraints?}
\newblock Mathematical Programming 131:37-48
(2012)

\bibitem{newoptcond}
S. Dempe, J. Dutta, and B.S. Mordukhovich,
\newblock {\it New necessary optimality conditions in optimistic bilevel programming},
\newblock Optimization 56 (5-6):577-604
(2007)

\bibitem{DempeZemkohoBook}
S. Dempe and A.B. Zemkoho (eds.), Bilevel optimization: advances and next challenges, {\small Springer (2020)}

\bibitem{bilevelreform}
S. Dempe  and A.B. Zemkoho,
\newblock {\it The bilevel programming problem: reformulations, constraint qualification and optimality conditions},
\newblock Mathematical Programming 138:447-473
(2013)

\bibitem{dempezemkoho1}
S. Dempe  and A.B. Zemkoho,
\newblock {\it The generalized Mangasarian-Fromowitz constraint qualification and optimality conditions for bilevel programs},
\newblock Journal of Optimization Theory and Applications 148(1):46-68 (2011)



\bibitem{FanYuan05}
J.Y. Fan and Y.X. Yuan,
\newblock {\it On the quadratic convergence of the Levenberg-Marquardt method without nonsingularity assumption},
\newblock Computing 74:23-39
(2005)

\bibitem{fischer92}
A. Fischer,
\newblock {\it A special Newton-type optimization method},
\newblock Optimization 24(3):269-284
(1992)


\bibitem{newtonbilevel18}
A. Fischer, A.B. Zemkoho, and S. Zhou,
\newblock {\it Semismooth Newton-type method for bilevel optimization: global convergence and extensive numerical experiments},
\newblock {arXiv:1912.07079}
(2019)

\bibitem{Fletcher75}
R. Fletcher,
\newblock {\it An ideal penalty function for constrained optimization},
\newblock IMA Journal of Applied Mathematics 15:319-342
(1975)


\bibitem{FTZ20}
J. Fliege, A. Tin, and A.B. Zemkoho,
\newblock {\it Gauss-Newton-type methods for bilevel optimization},
\newblock Computational Optimization and Applications, \href{https://doi.org/10.1007/s10589-020-00254-3}{doi.org/10.1007/s10589-020-00254-3}
(2021)





\bibitem{HerskovitsTanakaLeontiev2013}
J. Herskovits, M.T. Filho, and A. Leontiev,
{\it An interior point technique for solving bilevel programming problems},
 Optimization and Engineering 14(3):381-394 (2013)


\bibitem{kanzow1996}
C. Kanzow,
\newblock {\it Some noninterior continuation methods for linear complementarity problems},
\newblock SIAM Journal on Matrix Analysis and Applications 17(4):851-868
(1996)






\bibitem{LamparielloSagratellaNumerically2020} L. Lampariello and S. Sagratella,
\newblock {\it Numerically tractable optimistic bilevel problems}, Computational Optimization and Applications 76(2):277-303 (2020)

\bibitem{lin14}
G.-H. Lin, M. Xu, and J.J. Ye,
\newblock {\it On solving simple bilevel programs with a nonconvex lower level program},
\newblock Mathematical Programming 144(1-2):277-305 (2014)


\bibitem{MMZ20}
P. Mehlitz, L.I. Minchenko, and A.B. Zemkoho,
\newblock {\it A note on partial calmness for bilevel optimization problems with linear structures at the lower level},
Optimization letters, \href{https://doi.org/10.1007/s11590-020-01636-6}{doi.org/10.1007/s11590-020-01636-6}
(2020)

\bibitem{mitsos08}
A. Mitsos, P. Lemonidis, and P.I. Barton,
\newblock {\it Global solution of bilevel programs with a nonconvex inner program},
\newblock Journal of Global Optimization 42(4):475-513 (2008)

\bibitem{nocedal99}
J. Nocedal and S.J. Wright,
 Numerical optimization, Springer (1999)


\bibitem{paulavicius17}
R. Paulavicius, J. Gao, P. Kleniati, and C.S. Adjiman,
\newblock {\it BASBL: Branch-and-sandwich bilevel solver. Implementation and computational study with the BASBLib test sets},
{\footnotesize Computers \& Chemical Engineering 132 (2020)}

\bibitem{pillo89}
G. Di Pillo and L. Grippo,
\newblock {\it Exact penalty functions in constrained optimization},
\newblock SIAM journal on control and optimization 27 (6): 1333-1360, (1989)

\bibitem{PinedaByllingMorales2018}
S. Pineda, H. Bylling, and J.M. Morales,
{\it Efficiently solving linear bilevel programming problems using off-the-shelf optimization software},
Optimization and Engineering 19(1):187-211 (2018)


\bibitem{PTVF92}
W.H. Press, S.A. Teukolsky, W.T. Vetterling, and B.P. Flannery,
\newblock {\it Numerical recipes in C: the art of scientific computing (2nd Ed)},
\newblock Cambridge University Press (1992)





\bibitem{LMnumer20}
A. Tin and A.B. Zemkoho,
\newblock {\it Supplementary material for ``Levenberg-Marquardt method for bilevel optimization''},
\newblock School of Mathematical Sciences, University of Southampton, UK (2020)

\bibitem{wieseman13}
W. Wiesemann, A. Tsoukalas, P. Kleniati, and B. Rustem,
\newblock {\it Pessimistic bilevel optimization},
\newblock SIAM Journal on Optimization 23(1):353-380
(2013)

\bibitem{xu14}
M. Xu  and J.J. Ye,
\newblock {\it A smoothing augmented lagrangian method for solving simple bilevel programs},
\newblock Computational Optimization and Applications 59(1-2):353-377 (2014)


\bibitem{YF01}
N. Yamashita and M. Fukushima,
\newblock {\it On the rate of convergence of the Levenberg-Marquardt method},
\newblock Computing 15:237-249
(2001)

\bibitem{optcondbil95}
J.J. Ye  and D.L. Zhu,
\newblock {\it Optimality conditions for bilevel programming problems},
\newblock  Optimization 33:9-27
(1995)

\bibitem{ZemkohoZhouComparison2021}
A.B. Zemkoho and S. Zhou,
\newblock {\it Theoretical and numerical comparison of the Karush–Kuhn–Tucker and value function reformulations in bilevel optimization},
\newblock Computational Optimization and Applications
\newblock \href{https://doi.org/10.1007/s10589-020-00250-7}{doi.org/10.1007/s10589-020-00250-7}
(2021)

\bibitem{zemkohothesis}
A.B. Zemkoho,
\newblock {\it Bilevel programming: reformulations, regularity and stationarity},
\newblock PhD thesis, Department of Mathematics and Computer Science, TU Bergakademie Freiberg, Freiberg, Germany
(2012)

\bibitem{bolib17}
S. Zhou, A.B. Zemkoho, and A. Tin,
\newblock {\it BOLIB: Bilevel Optimization LIBrary of test problems},
\newblock In: S. Dempe and A.B. Zemkoho (eds.), Bilevel optimization: advances and next challenges,
\newblock Springer (2020)
\end{thebibliography}
\end{document}